\documentclass[a4paper,11pt]{article}
 \usepackage[latin1]{inputenc}   
\usepackage[T1]{fontenc}         
\usepackage{vmargin}
\setpapersize{A4}
\usepackage{amssymb,amsmath,amscd}

\def\A{\mathcal{A}}
\def\D{\mathcal{D}}
\def\F{\mathcal{F}}
\def\S{\mathcal{S}}
\def\T{\mathcal{T}}
\def\E{\mathbb{E}}
\def\N{\mathbb{N}}
\def\P{\mathbb{P}}
\def\Q{\mathbb{Q}}
\def\R{\mathbb{R}}
\def\Z{\mathbb{Z}}

\def\pentsup#1{\left\lceil #1 \right\rceil}
\def\un{\mathbf{1}}
\def\dd{\mathrm d}
\def\eps{\varepsilon}

\renewcommand{\le}{\leqslant}
\renewcommand{\ge}{\geqslant}

\newtheorem{lemma}{Lemma}
\newtheorem{definition}[lemma]{Definition}
\newtheorem{prop}[lemma]{Proposition}
\newtheorem{theorem}{Theorem}

\newtheorem{remark}[lemma]{Remark}
\newtheorem{corollary}[lemma]{Corollary}

\newenvironment{proof}[1][]%
{\noindent{\bf Proof{}#1:}\par\nobreak}{\nopagebreak\hspace*{\fill} $\square$\vspace*{1.4ex}\par}

\def\II{\mathbf I} 
\def\ab{\{a,  \ldots, b \}}
\def\lef{\mathsf l} \def\rig{\mathsf r}
\def\voi#1{\mathrm{Adj}(#1)}
\def\bxx{\mathbf x}
\def\prr{\rho} \def\pry{\eta} \def\prp{\pi}
\def\prt{\Pi}
\def\fr{F}
\def\AU{\mathcal A_{U}} \def\AV{\mathcal A_{V}}
\def\AW{\mathcal A_{W}} \def\AR{\mathcal A_{R}} \def\AQ{\mathcal A_{Q}}
\def\AZ{\mathcal A_{0}} 
\def\UU{\mathcal U} \def\VV{\mathcal V}
\def\WW{\mathcal W} \def\RR{\mathcal R} \def\QQ{\mathcal Q}
\def\PP{\mathcal P}
\def\typeu{U} \def\typev{V} \def\typew{W} \def\typer{R} \def\typeq{Q}
\def\umat{\mathbb U} \def\rmat{\mathbb W} 
\def\vmat{\mathbb V} \def\fmat{\mathbb F}

\def\sv{v_0} \def\sw{w_0} \def\ss{\sigma}

\def\vara{\kappa} \def\varb{\nu} \def\bata{\alpha}
\def\cpgoe{\mathrm{CpGo/e}}
\def\tpaoe{\mathrm{TpAo/e}}

\begin{document}
\title{Solvable models of neighbor-dependent nucleotide substitution processes}
\date{October 10, 2005}
\author{Jean B\'erard, Jean-Baptiste Gou\'er\'e, Didier Piau}

\maketitle

\bibliographystyle{plain}

\begin{abstract}
We prove that a wide class of models of Markov neighbor-dependent substitution
processes on the integer line is solvable. This class contains
some models of nucleotide substitutions recently introduced and studied
empirically by molecular biologists. We show that the polynucleotide
frequencies at equilibrium solve explicit finite-size linear
systems.  Finally, the dynamics of the process and the
distribution at equilibrium exhibit some
stringent, rather unexpected,
independence properties. For example, nucleotide 
sites at distance at least three evolve independently, and
the sites, if encoded as purines and pyrimidines,
evolve independently.
\end{abstract}

\section*{Introduction}
\addcontentsline{toc}{section}{Introduction}
\label{s.intro}
In simple Markov models of  nucleotide substitution
  processes, one assumes that each site along the DNA sequence evolves
  independently of the other sites and according to some specified
  rates of substitution. 
We introduce some notations, which are quite 
common for the biologist and useful 
to describe conveniently 
these models, as well as, later on, the more sophisticated 
ones which are the subject of this paper.
 For these definitions and, later on
in the paper, other nomenclatures, see~\cite{nome}.
\begin{definition}\label{d.d0}
The nucleotide
  alphabet is
$$
\A:=\{A,T,C,G\}.
$$ These letters stand for Adenine, Thymine, Cytosine and Guanine,
respectively.  The nucleotides $A$ and $G$ are purines, often
abbreviated by the letter $R$, the nucleotides $T$ and $C$ are
pyrimidines, often abbreviated by the letter $Y$.  The canonical
projection $\prp$ on the purine/pyrimidine alphabet $\{R,Y\}$ is
defined by
$$
\prp(A):=R=:\prp(G),\qquad
\prp(T):=Y=:\prp(C).
$$
Substitutions of the form $R\to R$ and $Y\to Y$ are called
transitions, substitutions of the form
$R\to Y$ and $Y\to R$ are called transversions.
Finally, for every subsets $X$ and $Z$ of $\A$, it is customary to
write 
$XpZ$ for the collection of
dinucleotides in $X\times Z$.
 \end{definition}
For instance,
YpR dinucleotides are formed by a purine followed by a pyrimidine,
hence there are $4$ such dinucleotides, namely,
CpG, TpA, TpG, and CpA.

Experimental facts are that
transitions are, in many cases, 
more frequent than transversions (typical ratios are 3:1), and 
that substitutions to
$C$  and to $G$ occur at different rates than substitutions to
$A$ and to $T$, see Duret and Galtier~\cite{dg} and the references in
their paper.
Many studies take this into account, for instance, 
in Tamura's model, one assumes that the matrix of substitution rates is
\begin{equation}\label{e.tam}
\begin{array}{rc}
\begin{array}{r} \\ A\\ T\\ C\\ G \end{array}
&
\!\!\!\!\begin{array}{c}
\begin{array}{cccc}A\, & \,T\, & \,C\, & \,G\end{array}
\\
\left(
\begin{array}{cccc}
\cdot & v_{2} & v_{1} & \kappa v_{1}
\\
v_{2} & \cdot & \kappa v_{1} & v_{1}
\\
v_{2} & \kappa v_{2} & \cdot & v_{1}
\\
\kappa v_{2} & v_{2} & v_{1} & \cdot
\end{array}
\right),
\end{array}
\end{array}
\end{equation}
where $v_1$, $v_2$ and $\kappa$ are nonnegative real numbers.
This means, for instance, that each
nucleotide $A$ is replaced by $T$ at rate $v_{2}$, by $C$ at rate $v_{1}$, and
by $G$ at rate $\kappa v_{1}$.

Under this model and under related ones with no interaction between
the sites, the nucleotide attached to any given site converges in
distribution to the stationary measure of the Markov chain described
by the matrix of the rates and, at equilibrium, the sites are
independent. This last consequence is unfortunate in a biological
context, since the frequencies $\fr(x)$ of the nucleotides and the
frequencies $\fr(xy)$ of the dinucleotides, observed in actual
sequences, are often such that, for many nucleotides $x$ and $y$,
$$
\fr(xy)\neq \fr(x)\,\fr(y).
$$
In fact, it is well known that the nucleotides in the immediate
neighborhood of a site
can affect drastically the substitution rates at this site.
For instance, in the genomes of vertebrates, the increased substitutions of
cytosine by thymine and of guanine by adenine in CpG dinucleotides
are often quite noticeable (typical ratios are 10:1, when compared to the
rates of the other substitutions).  The chemical reasons of this
CpG-methylation-deamination  
process are also well known and one can guess that, at equilibrium,
the number of CpG is decreased while the
number of TpG and CpA is increased when one adds 
high rates of CpG substitutions to Tamura's model.

The need to incorporate these effects into
more realistic models of nucleotide substitutions seems widely
acknowledged. However, the exact consequences of the introduction of such
neighbor-dependent substitution processes (in the case above, $CG\to
CA$ and $CG\to TG$), while crucial for a quantitative assessment of
these models, remain virtually unknown, at least up to our knowledge,
on a theoretical ground.  To understand why, note that the
distribution of the value of the nucleotide at site $i$ at a given
time depends a priori on the values at previous times of the
dinucleotides at sites $i-1$ and $i$ (to account for the $CG\to CA$
substitutions), and at sites $i$ and $i+1$ (to account for the $CG\to
TG$ substitutions), whose joint distribution, in turn, depends a
priori on the values of some trinucleotides, and so on.

Duret and Galtier~\cite{dg} introduced and analyzed a model, which we call
Tamura + CpG, that adds to Tamura's rates of substitution the
availability of substitutions $CG\to CA$ and $CG\to TG$, both at the
additional rate $\varrho\ge0$.  (These authors used a parameter
$\kappa_{1}\ge\kappa$, related to our parameter $\varrho$ and to the
ratio $\theta:=v_{1}/(v_{1}+v_{2})$, by the simple equation
$\varrho=(1-\theta)\,(\kappa_{1}-\kappa)$.)  To evade the curse,
explained above, of recursive calls to the frequencies of longer and
longer words, Duret and Galtier used as approximate frequencies
$F(xyz)$ of the trinucleotides the values
$$
F(xyz)\approx F(xy)\,F(yz)/F(y).
$$ Interestingly, these approximations would be exact if the sequences
at equilibrium followed a Markov model (with respect to the space
index $i$). This is not the case but, relying on computer simulations,
Duret and Galtier studied the $G+C$ content at equilibrium and other
quantities of interest in a biological context, they compared the
simulated values to the values predicted by their truncated model, and
they showed that their approximations captured some features of the
behavior of the true model.  In particular, they highlighted that
these models had inherent, and previously unexpected, consequences on
the frequency of the TpA dinucleotide as well.  We mention that Arndt
and his coauthors considered similar models and their biological
implications, see Arndt, Burge and Hwa~\cite{abh}, Arndt~\cite{arn}
and Arndt and Hwa~\cite{ah} for instance.

In this paper we introduce 
a wide extension of the Tamura + CpG model of neighbor-dependent
substitution processes, which we call R/Y + YpR models,
and we show that these models are solvable.
More precisely, we prove that the 
frequencies of polynucleotides at equilibrium
solve explicit finite-size linear systems. 
Thus, the infinite regressions to the frequencies of
longer and longer polynucleotides 
 described above disappear, and one can 
compute analytically several quantities of interest related to these models. 
For instance, one can assess rigorously 
 the effect of neighbor-dependent substitutions.
As noted above, the very possibility of such a solution comes as a
surprise. Additionally, our analysis provides some stringent
independence properties of these models at equilibrium, which, in our
view, should make the biologist somewhat 
more reluctant to use them.

\paragraph{Acknowledgements}
We wish to thank the biologists who contributed to this
work, in particular Laurent Duret, Nicolas Galtier, Manolo Gouy and Jean
Lobry. When the need arose, they willingly provided facts, references, and 
their numerous insights about the subject of this study.


\section{Description of the models}
\label{s.idotm}
\subsection{R/Y models}
\label{ss.rym}
To explain the basis of our analysis, 
we first note a key property of Tamura's model, described by matrices
of type
(\ref{e.tam})
in the
introduction.
\begin{quote}
\textbf{Property (a)} One can complete the matrix of the substitution
rates by diagonal elements such that, on the one hand, the joint
distributions of the elapsed times before a substitution occurs and of
the nucleotide that this substitution yields coincide for the two
purines, and, on the other hand, the joint distribution of the elapsed
times before a substitution occurs and of the nucleotide that this
substitution yields coincide for the two pyrimidines.
\end{quote}
\begin{definition}
\label{d.rsury}
The R/Y models are the models of substitutions such that property (a) holds. 
\end{definition}
Property (a) involves a comparison of some coefficients of the $A$ and
$G$ lines of the matrix of substitution rates, 
and of its $T$ and $C$ lines, namely, the coefficients that correspond
to the transversions.
In the most general model such that property (a) holds, the matrix
of substitution rates  is
$$
\begin{array}{rc}
\begin{array}{r} \\ A\\ T\\ C\\ G \end{array}
&
\!\!\!\!\begin{array}{c}
\begin{array}{cccc}A\, & \,T\, & \,C\, & \,G\end{array}
\\
\left(
\begin{array}{cccc}
\cdot & v_{T} & v_{C} & *
\\
v_{A} & \cdot & * & v_{G}
\\
v_{A} & * & \cdot & v_{G}
\\
* & v_{T} & v_{C} & \cdot
\end{array}
\right),
\end{array}
\end{array}
$$
where the parameters $v_x$ are nonnegative.
Here $*$ signals some free coefficients, which correspond to the
transitions. 
The following characterization of the R/Y models is immediate.
\begin{prop}[R/Y substitution rates]\label{p.crym}
A substitution model is R/Y if and only if there exists nonnegative rates
$v_{x}$ and $w_{x}$ such that the matrix of substitution rates is
$$
\begin{array}{rc}
\begin{array}{c} \\ A\\ T\\ C\\ G \end{array}
&
\!\!\!\!\begin{array}{c}
\begin{array}{cccc}A\, & \,T\, & \,C\, & \,G\end{array}
\\
\left(
\begin{array}{cccc}
\cdot & v_{T} & v_{C} & w_{G}
\\
v_{A} & \cdot & w_{C} & v_{G}
\\
v_{A} & w_{T} & \cdot & v_{G}
\\
w_{A} & v_{T} & v_{C} & \cdot
\end{array}
\right).
\end{array}
\end{array}
$$ 
\end{prop}
This means, for instance, that each nucleotide $A$ is replaced by $T$ at rate $v_{T}$, by $C$
at rate $v_{C}$, and by $G$ at rate $w_{G}$.  
The rates $v_x$ and $w_x$ are indexed by the nucleotide $x$ 
that the
corresponding substitution produces.

The full matrix of an R/Y model, which shows that
property  (a) holds, is
$$
\begin{array}{rc}
\begin{array}{c} \\ A\\ T\\ C\\ G \end{array}
&
\!\!\!\!\begin{array}{c}
\begin{array}{cccc}A\, & \,T\, & \,C\, & \,G\end{array}
\\
\left(
\begin{array}{cccc}
w_{A} & v_{T} & v_{C} & w_{G}
\\
v_{A} & w_{T} & w_{C} & v_{G}
\\
v_{A} & w_{T} & w_{C} & v_{G}
\\
w_{A} & v_{T} & v_{C} & w_{G}
\end{array}
\right),
\end{array}
\end{array}
$$ 
where $v_x$ and $w_x$ are nonnegative rates.
We recall that the diagonal elements represent fictitious substitutions $x\to x$, 
which leave the whole process unchanged.

\subsection{Situation of the R/Y models}
\label{ss.sotrym}
Here are some relationships between the R/Y model of substitutions 
and other classical ones, see~\cite{whe}.
First, Tamura-Nei's model (usually abbreviated as TN93) 
corresponds to the restriction of R/Y such that
$$
w_A\,v_G=w_G\,v_A,
\qquad
w_T\,v_C=w_C\,v_T.
$$
Hence special cases of TN93 are also special cases of R/Y.
For instance, Felsenstein's model (F84)
corresponds to the restriction of TN93 such that
$$w_A+v_T+v_C+w_G=v_A+w_T+w_C+v_G.
$$
Hasegawa-Kishino-Yano's model (HKY85)
corresponds to the restriction of R/Y such that
 $w_x/v_x$ does not depend on $x$.
In other words, both F84 and HKY85 
are subcases of TN93,
which is a subcase of R/Y.

Kimura's model with two parameters (K2P, also known as K80) 
corresponds to the restriction of R/Y such that
$v_x$ and $w_x$ do not depend on $x$.
Jukes-Cantor's model (JC69) corresponds to the restriction of R/Y such that
$v_x=w_x$ and $v_x$ does not depend on $x$.
Tamura's model, as it appears in the aforementioned paper by
Duret and Galtier, is intermediate between HKY85 and K2P since it
corresponds to the restriction of R/Y such that $v_A=v_T$,
$v_C=v_G$, and
$w_x/v_x$ does not
depend on 
$x$.

We summarize this as a proposition.
\begin{prop}\label{p.mod}
The model
JC69 is a strict subcase
of K80, which is a strict subcase of both HKY85 and F84, which are both strict
subcases of TN93. All these, and Tamura's model in Duret and Galtier (2000),
are strict subcases of R/Y.
\end{prop}
Finally, we mention that the
general time-reversible model (GTR) is not comparable with R/Y, 
in other words some matrices of rates of substitutions are GTR but not
R/Y, and vice versa. 
The intersection of GTR and R/Y is TN93. 
As a consequence, the complement of TN93 in R/Y contains only
non-reversible models.

\subsection{YpR substitutions}
\label{ss.yprs}
We generalize the CpG mechanism of substitutions considered by 
Duret and Galtier (2000), adding specific rates of substitutions from
each YpR dinucleotide as follows.
\begin{itemize}
\item
Every dinucleotide $CG$ moves to $CA$ at rate $r^C_{A}$
and to $TG$ at rate $r^G_{T}$.
\item
Every dinucleotide $TA$ moves to $CA$ at
rate $r^{A}_{C}$ and to $TG$ at rate $r^{T}_{G}$. 
\item
Every dinucleotide $CA$ moves to $CG$ at rate $r^C_{G}$
and to $TA$ at rate $r^A_{T}$.
\item
Every dinucleotide $TG$ moves to $CG$ at
rate $r^{G}_{C}$ and to $TA$ at rate $r^{T}_{A}$. 
\end{itemize}
The rationale for these notations is as follows. The rates $r^y_{x}$
are indexed by the nucleotide $x$ produced by the
substitution, and by the nucleotide $y$ completing the YpR
dinucleotide $xy$ that the substitution yields, when $x$ is a
pyrimidine, or completing the YpR dinucleotide $yx$ that the
substitution yields, when $x$ is a purine. In other words, $y$ is the
nucleotide not affected by the substitution.
\begin{definition}\label{d.full}
R/Y + YpR models correspond to the superposition of rates of
substitutions of an R/Y model and of rates of substitution $r^y_{x}$
of dinucleotides YpR, as described above.
\end{definition}
To recover the model of Duret and Galtier, one should assume that
$$
r^C_{A}=r^G_{T}, \quad r^{A}_{C}=r^T_{G}=r^C_{G}=r^{A}_{T}=r^G_{C}=r^T_{A}=0.
$$ R/Y + YpR models use $16$ mutation rates, namely, $4$ parameters
$v_{x}$ for the transversions, $4$ parameters $w_{x}$ for the
transitions, and $8$ parameters $r^y_{x}$ for the mutations involving
YpR dinucleotides.  To multiply these $16$ parameters by the same
scalar changes the time scale, but not the evolution itself nor the
stationary distributions, hence one can consider that the class R/Y +
YpR has dimension $15$.  We allow for possibly negative values of the
rates $r^{y}_{x}$ and the model makes sense if the following
inequalities are satisfied:
$$
v_{x}\ge0,\quad
w_{x}\ge0,\quad
w_{x}+r^{y}_{x}\ge0.
$$
Finally, note that the rates $r^{y}_{x}$ describe transitions, and never transversions.

\subsection{Sets of sites}
\label{ss.sos}
DNA sequences are represented by sequences of letters of the alphabet
$\A$. 
Although real DNA sequences are obviously finite, their typical length
is rather large, hence one often considers them as 
doubly infinite sequences of letters, that is, as 
elements of the set $\A^{\Z}$. Then sites on the DNA sequence are identified 
to integer numbers in $\Z$. 
Our mathematical results are most conveniently expressed in this
setting. However, as we explain below, 
most of them are valid for suitable finite sets of
sites as well.


\section{Description of the results}
\label{s.idotr}
Roughly speaking, our main findings about R/Y + YpR models at
equilibrium 
are the following.
\begin{itemize}
\item One can compute the exact
value of the frequency of every polynucleotide by solving a finite-size linear system.
\item One can simulate exact samples of
sequences of nucleotides.
\item Some surprising
independence properties between sites hold.
\end{itemize}
Furthermore, the dynamics satisfies the following.
\begin{itemize}
\item The dynamics converges to the unique equilibrium, for every starting distribution.
\end{itemize}
To be more specific, we prove the following results. 
\begin{theorem}[Construction and general properties]\label{t.constgen}
There exists a unique Markov process $(X(s))_{s \geq 0}$ on $\A^{\Z}$
associated to the transition rates defined in section~\ref{s.idotm}. 
Under a generic non-degeneracy condition (ND), stated in
section~\ref{s.nad}, this process is ergodic, that is, it has a unique
stationary distribution $\mu$ and, for any initial distribution,
$X(s)$ converges to $\mu$ in distribution as $s$ goes to
infinity. Moreover, $\mu$ is invariant and ergodic with respect to the
translations in $\Z$.
\end{theorem}

From theorem~\ref{t.constgen}, equilibrium properties of the model are
well-defined, for instance the equilibrium frequency of polynucleotides.

\begin{theorem}[Dynamics]\label{t.dyn}
There exists an i.i.d.\ sequence of
marked Poisson point processes $(\xi_i)_i$ on the real line, 
indexed by the sites $i$, 
and a measurable function $\Phi$ with
values in
$\A$ such that, for every couple of times $s\le t$ and
every site
$i$, the value $X_{i}(t)$ of site $i$ at time $t$ is
$$
X_{i}(t)=\Phi(X_{i-1}(s),X_{i}(s),X_{i+1}(s);\xi_{i-1}\cap[s,t],
\xi_{i}\cap[s,t],\xi_{i+1}\cap[s,t]).
$$
\end{theorem}
\begin{theorem}[Statics]\label{t.sta}
Assume that the distribution of $(X_{i})_{i}$ is stationary.
There exists an i.i.d.\ sequence of marked
Poisson point processes $(\xi_{i})_{i}$ on the real halfline, indexed by the sites $i$,
and a measurable function $\Psi$ with values in $\A$ such that, for every site $i$,
$$
X_{i}=\Psi(\xi_{i-1},\xi_{i},\xi_{i+1}).
$$
\end{theorem}
Theorems~\ref{t.dyn} and \ref{t.sta} describe some structural properties,
at the basis of our subsequent results.
We begin with some direct consequences.
\begin{definition}
\label{d.voi}
For every subset $I$ of $\Z$, let $\voi{I}$ denote the set of integers
$i$ such that either $i-1$ or $i$ or $i+1$ belongs to $I$.
\end{definition}
\begin{prop}[Dynamics]\label{p.dy}
Consider  sequences indexed by $\II$.
The restriction of the dynamics to a subset of sites
$I$ does not depend on $\II$, as soon as $\II$ contains $\voi{I}$.
\end{prop}
Proposition~\ref{p.dy} shows for instance that the behavior of the
sites in $\{1,\ldots,n\}$ is the same, whether one considers that
these sites are embedded in a sequence indexed by $\Z$, or in a
sequence indexed by $\{0,1,\ldots,n+1\}$. If the sites are embedded in
$\{0,1,\ldots,n+1\}$ or in a larger finite set, this means that the
boundary conditions have no effect on the behavior of the sites in
$\{1,\ldots,n\}$, hence one can consider at will discrete intervals or
discrete circles. For instance, one can decide, either that the only
neighbor of $0$ is $1$ and the only neighbor of $n+1$ is $n$, or that
$0$ has neighbors $1$ and $n+1$ and that $n+1$ has neighbors $n$ and
$0$. This decision will modify the evolution at the sites $0$ and
$n+1$ but not at the sites in $\{1,\ldots,n\}$.  This remark concerns
all the results that we state later on in this section.

We come back to the consequences of theorems~\ref{t.dyn} and
\ref{t.sta}.
\begin{corollary}[Statics]
Assume that the distribution of $(X_{i})_{i}$ is stationary.
Fix some sets of sites $I_k$ such that the sets $\voi{I_k}$ are disjoint. 
Then the families $(X_i)_{i\in I_k}$ are independent from each other.
\end{corollary}
The condition in the corollary means that, for every $k\neq k'$,
$|i-i'|\ge3$ for every $i$ in $I_k$ and every $i'$ in $I_{k'}$.  For
instance, the sequence $(X_{3i})_{i\in\Z}$ is i.i.d.\ at equilibrium.

Our next results deal with what is probably the main concern of biologists in
relation to this model, namely, the actual computation of some equilibrium
frequencies.
\begin{theorem}\label{t.ration}
The equilibrium frequencies of polynucleotides solve explicit finite-size linear
systems
and can be expressed as rational functions of the parameters of the model. 
\end{theorem}
\begin{theorem}[Nucleotides and YpR dinucleotides]\label{t.sn}
The frequency of each nucleotide at equilibrium can be expressed
explicitly as an affine function of the equilibrium frequencies of the
YpR dinucleotides.  Furthermore, the equilibrium frequencies of the
YpR dinucleotides solve an explicit $4\times 4$ linear system.
\end{theorem}
We state theorem~\ref{t.sn} more precisely as theorem~\ref{t.slin} in
section~\ref{ss.computYpR}.
Our last result is a consequence of the inner properties of our basic
construction.
\begin{theorem}[R/Y sequences]
\label{t.ry}
Encode the sequence of nucleotides as an $R/Y$ sequence of purines
and pyrimidines.
If the sequence of nucleotides at equilibrium, then
one obtains an i.i.d.\ $R/Y$ sequence with
 weights $t_R$ and $t_Y$, and
$$t_{Y}:=\frac{v_{C}+v_{T}}{v_A+v_T+v_C+v_G},\qquad
t_{R}:=\frac{v_{A}+v_{G}}{v_A+v_T+v_C+v_G}.
$$
\end{theorem}
In particular, the values of $t_Y$ and $t_R$ do not depend on
the values of the YpR substitution rates. 
This fact reflects, once again, the specific property of R/Y + YpR
substitution models that is at the basis of our analysis, namely, 
that the global model is equivalent to
the superposition of the
double substitution processes
described by the rates $r^y_x$ on top of the simple substitution
processes described by the rates $v_x$ and $w_x$, see section~\ref{s.nad}.
This remark shows that 
one cannot compute analytically, at least along these principles, 
the stationary measure
of substitution models that are 
not in the R/Y + YpR class, and in fact, we suspect that one cannot compute
it at all.

\section{Description of the paper}
\label{s.potrotp}

Here is a moderately detailed description of the content and
organization of the rest of the paper.

Part~A is devoted to a rigorous construction of the processes
described informally in section~\ref{s.idotm}.  We stress that one
could rely on the general principles in Liggett~\cite[chapter~1]{Lig},
based on infinitesimal generators, to build Markov processes on the
space of finite or infinite nucleotide sequences that correspond to
the jump rates defined above. However, in the present case, a direct
construction of the dynamics is possible, which yields straightforward
proofs of important structural properties of the process and possesses
some interesting coupling properties.  In section~\ref{s.nad}, we give
some notations and definitions.  In section~\ref{s.fb}, we explain the
construction of the process when the number of sites is finite. In
section~\ref{s.tpoaz}, we deduce from this the case when the number of
sites is infinite.  Finally, section~\ref{s.tryp} deals with the
simplifications related to the encoding of the nucleotide sequences as
sequences of purines and pyrimidines.

Part~B is devoted to the actual computation of some quantities of
interest for the model at equilibrium. From structural properties of
the process, described in part~A, computing the equilibrium
frequencies of polynucleotides amounts to solving finite-size linear
systems.  Moreover, one can express the nucleotide frequencies as
functions of the YpR dinucleotide frequencies, and these, in turn,
solve a $4 \times 4$ linear system. Section~\ref{s.tgc} explains this
in the general case.  

The remaining sections of part~B are devoted to restricted versions of
the model, which involve a reduced number of free parameters.  This
avoids, first, cumbersome formulas that would depend on a large number
of free parameters, and second, the prohibitive computational burden
involved in symbolically solving large linear systems.
Section~\ref{s.ussr} deals with uniform simple substitution rates and
double substitution rates from CpG and TpA only.  Section~\ref{s.wms}
deals with models that are invariant by the classical symmetry of DNA
strands. This case has biological significance since the invariance
reflects the well known complementarity of the two strands of DNA
molecules.  Section~\ref{s.soa} deals with the simplest non trivial
version of the R/Y + YpR model, namely, the case when all the simple
substitution rates coincide and there are no double substitutions
except from CpG to CpA and TpG, both at the same rate.  In this
setting, we provide the exact $16$ dinucleotide frequencies (that is,
not only the $4$ YpR frequencies), a perturbative analysis of the
frequency of every polynucleotide at vanishing CpG substitution rates,
and we describe the non degenerate limit of the model at high CpG
substitution rates.
The relatively shorter section~\ref{s.cd} 
explains the dynamics of the model, that is, its
evolution to the stationary measure.
The results of this section
are valid in our general setting but we expose them in the simplest
case.

Part~C deals with the actual simulation of these systems.
Coupling-from-the-past techniques are pivotal here, as explained in
section~\ref{s:CFTP}.  Section~\ref{s:bepractical} delves into the
details of the effective implementation of the CFTP method in the
context of R/Y + YpR systems, and provides the basic schemes of two
different algorithms that perform this sampling method. Once again,
the crucial fact here is that one can use finite-size sets of sites to
simulate the behavior of smaller finite-size sets of sites.

Finally, section~\ref{s.ion} provides an index of the main parameters
and notations used in the paper.

\section*{Part A Construction}
\addcontentsline{toc}{section}{Part A Construction}

\section{Notations and definitions}
\label{s.nad}

For every topological space $E$ and every real number $\sigma$,
$\D(\sigma, E)$ denotes the space of c\`adl\`ag (right continuous with
left limits) functions from $[\sigma,+\infty)$ to $E$, equipped with
the Skorohod topology and the corresponding Borel $\sigma$-algebra.

Let $\II$ denote the collection of nucleotide sites, thus $\II$ may be
either the integer line $\Z$ or a finite interval of integers.
Due
to technical reasons that will become apparent when we explain the
construction of the dynamics, $\II$ must contain at least $3$ sites.

We recall from definition~\ref{d.d0} in the introduction that $\A :=
\{ A,T,C,G\}$ denotes the nucleotidic alphabet, that $A$ and $G$ are
purines, encoded by $R$, that $C$ and $T$ are pyrimidines, encoded by
$Y$, and that the mapping $\prp$ is defined on $\A$ by
$$
\prp(A):=R=:\prp(G),\quad\prp(C):=Y=:\prp(T).
$$
We now define other mappings on $\A$.
\begin{definition}
\label{d.d1}
Let $\prr$ denote the application which fuses the two
purines together, and $\pry$ the application which fuses
the two pyrimidines
together, that is
$$
\prr(A):=R=:=\prr(G), \quad \prr(C):=C,\quad\prr(T):=T,
$$
and
$$
\pry(A):=A,\quad\pry(G):=G,\quad\pry(C):=Y=:\pry(T).
$$ 
For every nucleotide $x$ in $\A$, let $x^\ast$ denote the unique
nucleotide such that $\{x,x^\ast\}=\{A,G\}$ or
$\{x,x^\ast\}=\{C,T\}$. In other words,
$x\mapsto x^\ast$ is the involution of $\A$ such that
$$
A^\ast:=G,\quad
T^\ast:=C,\quad
C^\ast:=T,\quad
G^\ast:=A.
$$
\end{definition}
For every subsets $K$ and $J$ of the integer line such that $J\subset
K$, let $\prt^{J}_{K}$ denote the canonical projection from $\A^{K}$
to $\A^{J}$.  When $J=\ab$ for two integers $a\le b$, we often omit
the mention of $K$ and write $\prt^{a,b}$ for $\prt^{J}_{K}$. In other
words, if $K$ contains $\ab$ and if $\bxx:=(x_k)_{k\in K}$, then
$$
\prt^{a,b}(\bxx):=(x_k)_{a\le k\le b}.
$$ 
Let $\S^{+}$ denote the set of countably infinite locally finite
subsets of the real line $\R$, and let $\S:=\S^{+}
\cup\{\emptyset\}$. We equip $\S$ with the usual $\sigma$--algebra in
the context of point processes, namely, the smallest $\sigma$-algebra
such that, for every Borel subset $A$ of the real line, the function
$N \mapsto \mbox{card} (N \cap A)$ is measurable.

We consider collections $\xi$ of elements of $\S$, defined as follows.
Let  $\AZ$ denote the disjoint union of five copies of $\A$, say 
$$\AZ:=\AU\cup\AV\cup\AW\cup\AR\cup\AQ.
$$
Let $\Omega_{0}$ denote the space 
$$
\Omega_{0}:=\S^{\AZ\times\II}.
$$ 
We call  $\F_{0}$ the product $\sigma$--algebra on $\Omega_{0}$ inherited from that of $\S$.
 Let $\xi$ be a generic element of $\Omega_{0}$.
Hence $\xi$ may be written as 
$$
\xi=:(\xi_{i})_{i\in \II},\quad\mbox{where each}\
\xi_{i}\ \mbox{belongs to}\ \S^{\AZ}.
$$ 
For every $x$ in
$\A$, let $\UU^{x}_{i}(\xi)$ denote the $x$--coordinate of $\xi_{i}$ in $\AU$, hence
$\UU^{x}_{i}(\xi)$ belongs to $\S$.
Likewise, $\VV^{x}_{i}(\xi)$, $\WW^{x}_{i}(\xi)$, $\RR^{x}_{i}(\xi)$, 
and $\QQ^{x}_{i}(\xi)$ respectively denote the $x$--coordinates of
$\xi_{i}$ in $\AV$, $\AW$, $\AR$, and $\AQ$, and belong to $\S$ as well. 
In other words, for every $\xi$ in $\Omega_{0}$ and
every $i$ in $\II$,
$$
\xi_{i}=:\big(\UU^{x}_{i}(\xi),\VV^{x}_{i}(\xi),\WW^{x}_{i}(\xi),
\RR^{x}_{i}(\xi),\QQ^{x}_{i}(\xi)\big)_{x\in\A}.
$$
Recall that, for every real number $r$, the positive part $r^{+}$ and the negative part $r^{-}$ 
of $r$ are both nonnegative and defined by
the
 relations $r=r^{+}-r^{-}$ and $|r|=r^{+}+r^{-}$. 
\begin{definition}
\label{d.crs}
For every nucleotide $x$ in $\A$, the combined rate of substitution $c_{x}$ is the
nonnegative real number defined as
$$
c_{x}:=w_{x}-\max\{(r^{y}_{x})^{-},(r^{y^\ast}_{x})^{-}\},
$$
where $\{y,y^\ast\}=\{C,T\}$ if $\pi(x)=R$
and  $\{y,y^\ast\}=\{A,G\}$ if $\pi(x)=Y$.
\end{definition}
Finally, the probability measure $\Q$ on $(\Omega_{0},
\F_{0})$ is such that, for every site $i$ in $\II$ and every nucleotide $x$ in $\A$,
\begin{itemize}
\item $\UU^{x}_{i}$ is a  homogeneous Poisson process on the real line with constant rate $\min(v_{x},   
c_{x})$,
\item $\VV^{x}_{i}$ is a homogeneous Poisson process on the real line with constant rate $(v_{x}-c_{x})^{+}$,
\item $\WW^{x}_{i}$ is a homogeneous Poisson process on the real line with constant rate $(c_{x}-v_{x})^{+}$,
\item $\RR^{x}_{i}$ is a homogeneous Poisson process on the real line with constant rate 
$|r^{y}_{x}|$, where the rate $r^{y}_{x}$ corresponds to YpR substitutions starting from CpG or TpA,
\item $\QQ^{x}_{i}$ is a homogeneous Poisson process on the real line with constant rate 
$|r^{y}_{x}|$, where the rate $r^{y}_{x}$ corresponds to YpR substitutions starting from TpG or CpA.
\end{itemize}
Thus, $\Q$ is uniquely specified by the following additional condition.
\begin{itemize}
\item The Poisson processes $\UU^{x}_{i}$, $\VV^{x}_{i}$, $\WW^{x}_{i}$, $\RR^{x}_{i}$, 
and $\QQ^{x}_{i}$,  for every site $i$ and every
nucleotide $x$, are independent.
\end{itemize}
One sees that every $\UU^{x}_{i}\cup\VV^{x}_{i}$ is a homogeneous Poisson process 
with constant rate $v_{x}$ and that every $\UU^{x}_{i}\cup\WW^{x}_{i}$ is a homogeneous Poisson process 
with constant rate $c_{x}$.

We now provide a brief and intuitive description of the construction
of the dynamics of the process, using these Poisson processes.  As is
usual, the points in the processes $\xi_{i}$ are the ringing times of
exponential clocks that rule the evolution of the sites. There exists
five types of moves, labelled as \typeu, \typev, \typew, \typer, and
\typeq.
\begin{itemize}
\item
Type \typeu.  When an exponential clock attached to $\UU^{x}_{i}$
rings, the nucleotide at site $i$ moves unconditionally to the value
$x$.
\item Type \typev.  When an exponential clock attached to
$\VV^{x}_{i}$ rings, the nucleotide at site $i$ moves to the value $x$
provided that this move corresponds to a transversion.
\item
Type \typew.  When an exponential clock attached to $\WW^{x}_{i}$
rings, the nucleotide at site $i$ moves to the value $x$ provided that
this move corresponds to a transition.
\item
Type \typer.  When an exponential clock attached to $\RR^{x}_{i}$
rings, the nucleotide at site $i$ moves to the value $x$ in the
following cases: if this move
corresponds to a YpR substitution from $CG$ or $TA$ when the
associated rate $r^{y}_{x} \ge 0$, and if this move corresponds to a
transition but not to a substitution from $CG$ or $TA$ when
$r^{y}_{x}<0$.
\item
Type \typeq.  When an exponential clock attached to $\QQ^{x}_{i}$
rings, the nucleotide at site $i$ moves to the value $x$ in the
following cases: if this move corresponds to a YpR substitution from
$TG$ or $CA$ when the associated rate $r^{y}_{x} \ge 0$, and if this
move corresponds to a transition but not to a substitution from $TG$
or $CA$ when $r^{y}_{x}<0$.
\end{itemize}
The rates of the Poisson processes are chosen in order to couple as
strongly as possible the transitions and the transversions that yield
the same nucleotide and to take properly into account the inhibitory
effect of YpR mutations when some rates $r^{y}_{x}$ are negative.

We introduce a subset $\Omega_{1}$ of $\Omega_{0}$, defined by the following conditions.
\begin{itemize}
\item The sets $\PP^{x}_{i}$ are disjoint, for every nucleotide $x$, every site
$i$, and for every symbol $\PP$ in the set  $\{\UU,\VV,\WW,\RR,\QQ\}$.
\item For every site $i$, there exists a symbol $\PP$ in the set  $\{\UU,\VV,\WW,\RR,\QQ\}$ and 
a nucleotide $x$ in $\A$ such that the set $\PP^{x}_{i}$ is infinite. 
\end{itemize}

We also introduce a non-degeneracy condition.
\begin{itemize}
\item
(ND)
For every nucleotide $x$ in $\A$, $v_{x}$ and $c_{x}$ are positive.
\end{itemize}
Under condition (ND), $\Q(\Omega_{0} \setminus \Omega_{1})=0$. 
To avoid the handling of tedious
exceptions, we assume from now on that condition (ND) holds
and we work exclusively on $\Omega_{1}$,
equipped with the Borel $\sigma$--field and the probability measure 
that are induced by those of $(\Omega_{0},\F_{0},\Q)$. 
We denote this new probability space 
by $(\Omega_{1},\F_{1},\Q)$.

\section{Construction on finite intervals}
\label{s.fb}

Before considering the full integer line, we define the dynamics on
finite discrete segments, with periodic boundary conditions.  The
choice of boundary conditions is somewhat arbitrary, and one could use
instead free boundary conditions or fixed boundary
conditions. However, the dynamics with periodic boundary conditions is
invariant by the translations of the discrete circle, and this fact
turns out to be quite useful since it reduces the dimension of the
linear system which yields the invariant distribution, see
section~\ref{s.tgc}.

\subsection{Definitions and notations}

In this whole section~\ref{s.fb}, we fix two integers $a$ and $b$ such
that $a+2\le b$ and we assume that
$$
\II:=\ab.
$$
The definitions below depend on the choice of $\II$ but, to alleviate the notations, we do not
always mention explicitly the dependence.

For every site $i$ in $\II$, we introduce  $\lef(i)$ as the neighbor of $i$ to the left of $i$ and $\rig(i)$ as the neighbor
of $i$ to the right of $i$. More precisely,
$$\lef(i) := i-1\ \mbox{if}\ i \neq a,\quad\lef(a) := b,
\quad\rig(i) := i+1\ \mbox{if}\ i
\neq b,\quad\rig(i) := a. 
$$ 
Since $a+2\le b$, for every site $i$,
the sites $i$, $\lef(i)$ and $\rig(i)$ are three different sites.

Fix a sequence $\bxx=(x_{i})_{i \in \II}$ in $\A^{\II}$ and a site $i$
in $\II$.  Our next definitions are related to the moves of types
\typeq\ and \typer, defined in section~\ref{s.nad}.  Say that:
\begin{itemize}
\item The sequence $\bxx$ accepts the substitutions of type \typer\ to $G$ at site $i$ if $x_{\lef(i)}x_{i} =TA$.
\item The sequence $\bxx$ accepts the substitutions of type \typer\ to $C$ at site $i$ if $x_{i}x_{\rig(i)} = TA$.
\item The sequence $\bxx$ accepts the substitutions of type \typer\ to $A$ at site $i$ if $x_{\lef(i)}x_{i} = CG$.
\item The sequence $\bxx$ accepts the substitutions of type \typer\ to $T$ at site $i$ if $x_{i}x_{\rig(i)} = CG$.
\end{itemize}
Likewise:
\begin{itemize}
\item The sequence $\bxx$ accepts the substitutions of type \typeq\ to $G$ at site $i$ if $x_{\lef(i)}x_{i} =CA$.
\item The sequence $\bxx$ accepts the substitutions of type \typeq\ to $C$ at site $i$ if $x_{i}x_{\rig(i)} = TG$.
\item The sequence $\bxx$ accepts the substitutions of type \typeq\ to $A$ at site $i$ if $x_{\lef(i)}x_{i} = TG$.
\item The sequence $\bxx$ accepts the substitutions of type \typeq\ to $T$ at site $i$ if $x_{i}x_{\rig(i)} = CA$.
\end{itemize}
Although this terminology makes little concrete sense when some
rates $r^{y}_{x}$ are negative, we use it even then.

\subsection{Construction}
\label{ss.cotp}
The goal of this section is to define, for  every real number
$\sigma$,
a measurable map
$$\varphi^{\II}_{\sigma} \ : \ \A^{\II} \times \Omega_{1}  \to  \D(\sigma, \A^{\II}).$$
Intuitively, each function $\varphi^{\II}_{\sigma}(\bxx,\xi)$
describes the dynamics 
in $ \A^{\II}$ that starts from the configuration $x$ at time
$\sigma$ and uses the moves prescribed by the realization
$\xi$ of the Poisson processes.
To be specific about the notations,
we write $$\varphi^{\II}_{\sigma}(\bxx,\xi)(s) = \left( \varphi^{\II}_{\sigma}(\bxx,\xi,i,s) \right)_{i \in \II}.
$$ 
In other words,
$\varphi^{\II}_{\sigma}(\bxx,\xi,i,s)$ stands for the $i$th coordinate
of the value of the function $\varphi^{\II}_{\sigma}(\bxx,\xi)$ at time $s\ge\sigma$, 
hence $\varphi^{\II}_{\sigma}(\bxx,\xi,i,s)$ belongs to $\A$.

From now on, we fix $\xi$ in $\Omega_{1}$ and $\bxx$ in $\A^{\II}$,
and, to alleviate the notations, we omit to mention the dependence with respect to $a$, $b$, $\sigma$ or $\xi$ of
various quantities. Introduce
$$\T:=\bigcup_{i\in\II}\T(i),
\quad
\T(i) :=  (\sigma,+\infty)\cap\bigcup_{z\in\A}
\UU_{i}^{z}\cup\VV^{z}_{i}\cup\WW^{z}_{i}\cup\RR_{i}^{z}\cup\QQ_{i}^{z}.
$$ Let $t_{-1}:=\sigma$ and $(t_{0} , t_{1},\cdots)$ denote the
ordered list of points in $\T$, that is,
 $$ \sigma < t_{0} < t_{1} < \cdots ,\quad\mbox{and}\quad \T =:
\{t_0,t_1,t_2,\ldots\}.
$$
For every $n \ge  0$, let $c_{n}$ denote the site where the $n$th move
occurs that affects a site in $\II$ after the time
$\sigma$, and let $M_{n}$ denote the description of this move. That is,
$$c_{n}=i\ \mbox{and}\
M_{n}=(z,\typeu)\ \mbox{if}\
t_{n}\ \mbox{belongs to}\ \UU_{i}^{z}.
$$ Likewise, $c_{n}=i$ and $M_{n}=(z,\typev)$, $M_{n}=(z,\typew)$,
$M_{n}=(z,\typer)$, and $M_{n}=(z,\typeq)$ respectively, if $t_{n}$
belongs to $\VV_{i}^{z}$, $\WW_{i}^{z}$, $\RR_{i}^{z}$, and
$\QQ_{i}^{z}$ respectively.  We often consider that $M_{n}$ belongs to
$\AZ$, for instance $M_{n}=(z,\typeu)$ may be identified with $z$ in
$\AU$. Alternatively, the second component of $M_{n}$ is considered as
a flag, it takes its values in the set
$\{\typeu,\typev,\typew,\typer,\typeq\}$, and it is often denoted by $f$.
In any case, the variables $c_{n}$ and $M_{n}$ are uniquely defined
for every $\xi$ in $\Omega_{1}$.
 
 We now define  a map 
$$\gamma_{\II} \ : \ \A^{\II} \times  \AZ \times \II \to \A.
$$
Fix  a sequence $\bxx:=(x_{i})_{i\in\II}$ in $\A^{\II}$,  a move description $M=(z,f)$ in $\AZ$, and a site $i$ in $\II$.
The definition of the $i$--coordinate of
$\gamma_{\II}(\bxx,z,f,i)$ differs from the definition of the other coordinates. More precisely,
one sets
 $$
\gamma_{\II}(\bxx,z,f,i)_{i} := z,
$$ 
if one of the following conditions is met.
\begin{itemize}
\item The flag $f =\typeu$.
\item The flag  $f =\typev$, and $x_{i}$ is a purine and $z$ is a pyrimidine,
  or vice versa. 
\item The flag $f =\typew$, and  $x_{i}$ and $z$ are both purines or both pyrimidines, 
\item The flag $f=\typer$ and the type \typer\ rate $r^{y}_{z} > 0$ and $\bxx$ accepts the substitutions of type \typer\
 to $z$ at site $i$. 
\item The flag $f=\typer$ and the type \typer\ rate $r^{y}_{z} < 0 $ and $\bxx$ does not accept the substitutions 
of type \typer\ to $z$ at site $i$ 
and $x_{i}$ and $z$ are both purines or both pyrimidines.
\item The flag $f=\typeq$ and the type \typeq\ rate $r^{y}_{z} > 0$ and $\bxx$ accepts the substitutions of type \typeq\
 to $z$ at site $i$.
\item The flag $f=\typeq$ and the type \typeq\ rate $r^{y}_{z} < 0 $ and $\bxx$ does not accept the substitutions 
of type \typeq\ to $z$ at site $i$ 
and $x_{i}$ and $z$ are both purines or both pyrimidines.
\end{itemize}
In every other case, that is, if $j=i$ and none of the conditions above is met, or if $j\neq i$, 
one sets
$$
\gamma_{\II}(\bxx,z,f,c)_{j}   :=   x_{j}.
$$ 
This defines the map $\gamma_{\II}$.
We now construct the map $\varphi^{\II}_{\sigma}$, using an  induction over increasing values of
the time $s\ge\sigma$. 
The initial condition is that, for  every $\sigma\le s< t_{0}$,
$$
\varphi^{\II}_{\sigma}(\bxx,\xi)(s) := \bxx.
$$ 
Assume now that  $\varphi^{\II}_{\sigma}(\bxx,\xi)(s)$ is well-defined
for every time $s$ such that $\sigma\le s<t_{n}$. 
Then, for every time $s$ such that $t_{n}\le s<t_{n+1}$, one sets
$$ 
\varphi^{\II}_{\sigma}(\bxx,\xi)(s) := 
\gamma_{\II}\left(\varphi^{\II}_{\sigma}(\bxx,\xi)(t_{n-1}),
M_{n},c_{n}\right).
$$
This defines a configuration 
$\varphi^{\II}_{\sigma}(\bxx,\xi)(s)$ for every time $s\ge\sigma$.
Since $\II$ is finite, measurability issues are obvious here.

\subsection{First properties}
\label{ss.prop}
Recall that in this section~\ref{s.fb}, $\II$ denotes the finite interval $\II=\ab$.
We record some immediate properties of the family of maps $\varphi^{\II}_{\sigma}$,
making use of still another notation.
\begin{definition}\label{d.xit}
For
every time set $T$, let
$
\xi T:=(\xi_{i}T)_{i\in\II}$,
where
$$
\xi_{i}T:=\big(\UU_{i}^{z}(\xi)\cap T,\VV_{i}^{z}(\xi)\cap T,\WW_{i}^{z}(\xi)\cap T,
\RR_{i}^{z}(\xi)\cap T,\QQ_{i}^{z}(\xi)\cap T\big)_{z\in\A}.
$$
\end{definition}
\begin{prop}\label{p.pf}
{\bf (1)} For every $s\ge t>\sigma$,
$$\varphi^{\II}_{\sigma} (\bxx, \xi)(s) =  \varphi^{\II}_{t}\left(   
\varphi^{\II}_{\sigma} (\bxx, \xi)(t) , \xi \right)(s).
$$
{\bf (2)}
For every  $s \ge  0$ and every real number $t$,
$$
\varphi^{\II}_{\sigma+t} (\bxx, \xi+t)(\sigma+t+s) = \varphi^{\II}_{\sigma} (\bxx, \xi)(\sigma+s),
$$
where $\xi+t$ denotes the result of the
addition of $t$ to each component of $\xi$. 
\\
{\bf (3)}
Finally, $\varphi^{\II}_{\sigma} (\bxx, \xi)(s)$ depends on $\xi$ only through $\xi[\sigma,s]$.
\end{prop}
The function $\varphi^{\II}_{0}(\bxx,\cdot)$, defined by
 $$
s \mapsto
\varphi^{\II}_{0}(\bxx, \cdot) (s),
$$
is a random variable on $(\Omega_{1},\F_{1},\Q)$ with values in $\D(0, \A^{\II})$.
Let $\P^{\II}_{\bxx}$ denote its distribution.
Then, from proposition~\ref{p.pf} and from the translation invariance of Poisson processes, 
 the family of measures
$$
\left\{\,\P^{\II}_{\bxx},\,\bxx \in \A^{\II}\,\right\}
$$ 
defines a Markov process in the sense of Liggett~\cite[chapter~1]{Lig}.
Moreover, it corresponds to the specification of the process on $\A^{\II}$ given in 
our section~\ref{s.idotm} in terms of transition rates, with periodic
boundary conditions. 

From now on, we use the notation 
$$
X^{\II}_{\bxx}(s):= \varphi^{\II}_{0}(\bxx, \cdot) (s).
$$ 
\begin{prop}\label{p.ef}
For every finite interval  $\II$ and every initial configuration $\bxx$, 
the Markov process $(X^{\II}_{\bxx}(s))_{s \ge  0}$ is ergodic. In other words, there exists a
unique stationary distribution $\mu_{\II}$ on $\A^{\II}$, and
for every $\bxx$ in $\A^{\II}$, $X^{\II}_{\bxx}(s)$ converges in distribution to
$\mu_{\II}$ when $s\to +\infty$.
\end{prop}

\begin{proof}[\ of proposition~\ref{p.ef}]
Immediate since $(X^{\II}_{\bxx}(s))_{s \ge  0}$ lives on a finite state space and is irreducible 
from the non-degeneracy
assumption (ND) at the end of section~\ref{s.nad}.
\end{proof}

\subsection{Dependencies}
\label{ss.depe}
We now make a fundamental observation.
\begin{prop}\label{p.i}
For every initial sequence $\bxx:=(x_{i})_{i\in\II}$, every
site $i$ in $\II$ and every time $s\ge\sigma$, the nucleotide
$\varphi^{\II}_{\sigma} (\bxx,
\xi,i,s)$, which depends a priori on all the information contained
in $\bxx$ and $\xi$, is in fact measurable with respect to the following restricted 
initial conditions and restricted sources of moves:
$$
x_{\lef(i)},\quad
x_{i},\quad
x_{\rig(i)},\quad
\xi_{\lef(i)}[\sigma,s],\quad
\xi_{i}[\sigma,s],\quad
\xi_{\rig(i)}[\sigma,s].
$$
More precisely,
there exists measurable maps $\Theta_{\sigma,s}$, independent of $a$ and $b$, and such that, for every integers $a$ and
$b$ such that $a+2
\le b$, and for every $a\le i \le b$ and $s \ge \sigma$,
$$
\varphi^{a,b}_{\sigma} (\bxx, \xi,i,s)    =    
\Theta_{\sigma,s}(x_{\lef(i)}, x_{i}, x_{\rig(i)},\xi_{\lef(i)}[\sigma,s], 
\xi_{i}[\sigma,s] ,
\xi_{\rig(i)}[\sigma,s]).
$$
\end{prop}
A straightforward consequence of proposition~\ref{p.i} is
proposition~\ref{p.r} below.
\begin{prop}\label{p.r}
Let $a$, $b$, $c_{1}$, $d_{1}$, $c_{2}$ and $d_{2}$ denote integers such that
$$c_k+1\le a\le b\le d_k-1, \quad
 k=1,2.
$$
Fix some initial conditions  $\bxx_{1}$ in $\A^{\{c_{1},\ldots,d_{1}\}}$
and $\bxx_{2}$ in $\A^{\{c_{2},\ldots,d_{2}\}}$ that
coincide on the interval $\{a-1,a,\ldots,b,b+1\}$,
that is,
such that
$$
\prt^{a-1,b+1}(\bxx_{1})    
=\prt^{a-1,b+1}(\bxx_{2}).
$$ 
Then, for every $s \ge  \sigma$,
$$  
\prt^{a,b} \left(\varphi^{c_{1},d_{1}}_{\sigma} (\bxx_{1}, \xi)(s)\right) = \prt^{a,b}  \left(\varphi^{c_{2},d_{2}}_{\sigma}
(\bxx_{2}, \xi)(s)\right).
$$
\end{prop}

\subsection{Proofs}
The proof of proposition~\ref{p.i} relies on lemma~\ref{l.d1} below. 
Recall definition~\ref{d.d1} of $\prr$ and $\pry$.
\begin{lemma}\label{l.d1}
For every site $i$ in $\II$ and every time $s \ge  \sigma$, 
the functions 
$$\varrho\left[\varphi^{\II}_{\sigma} (\bxx,
\xi,\lef(i),s)\right],\quad
\varphi^{\II}_{\sigma} (\bxx, \xi,i,s),\quad
\eta\left[\varphi^{\II}_{\sigma} (\bxx,
\xi,\rig(i),s)\right],
$$
which depend a priori on all the information in $\bxx$ and $\xi$, are in fact measurable with respect to
the following restricted initial conditions and restricted source of moves:
$$
\varrho(x_{\lef(i)}),\quad
x_{i},\quad
\eta(x_{\rig(i)}),\quad
\xi_{\lef(i)}[\sigma,s],\quad
\xi_{i}[\sigma,s],\quad
\xi_{\rig(i)}[\sigma,s].
$$ 
More precisely,
one may define a measurable map $\Psi_{\sigma,s}$, independent of $a$ and $b$, and such that, 
for every integers $a$ and
$b$ such that $a+2\le b$, every $a\le i\le b$, 
and every $s\ge\sigma$,
the triple
$$ 
\left( \varrho\left[\varphi^{a,b}_{\sigma} (\bxx, \xi,\lef(i),s)\right] ,  
 \varphi^{a,b}_{\sigma} (\bxx, \xi,i,s) , 
\eta\left[\varphi^{a,b}_{\sigma} (\bxx,\xi,\rig(i),s)\right]
\right)
$$
coincides with
$$
\Psi_{\sigma,s} \left( \varrho(x_{\lef(i)}), x_{i}, \eta(x_{\rig(i)}) ,\xi_{\lef(i)}[\sigma,s],  \xi_{i}[\sigma,s] ,
\xi_{\rig(i)}[\sigma,s] \right).
$$
\end{lemma}
\begin{proof}[\ of lemma~\ref{l.d1}]
Fix a source of moves $\xi$ in $\Omega_{1}$, an initial configuration $x$ in $\A^{\II}$,
a site $i$ in $\II$, and  let
$$
\T^* =  \T(\lef(i))\cup\T(i)\cup\T(\rig(i)),
$$
where the sets $\T(\cdot)$ are defined in section~\ref{ss.cotp}.
Let  $(t^*_{0} , t^*_{1},\cdots)$ denote the ordered list of
points in $\T^*$, that is,
 $$
\T^* =
\{t^*_0,t^*_1,t^*_2,\ldots\}\quad\mbox{where}\quad\sigma < t^*_{0} < t^*_{1} < \cdots,
$$
 and set $t^*_{-1}=\sigma$.
Then, for all $n \ge  0$, we describe the move that occurs at time
 $t^*_{n}$ through the description $M^*_{n}=(z^*_{n},f^*_{n})$ of this
 move
and the site $c^*_{n}$ where the move 
occurs, as defined in section~\ref{ss.cotp}.
Note that the moves that affect the sites $\lef(i)$, $i$ and  $\rig(i)$, 
 can only occur at one of the times $t^*_{n}$ for $n\ge0$.
 
Thus, to prove lemma~\ref{l.d1}, it is enough to prove that, for all $n\ge  -1$, 
 $$
\varrho\left[\varphi^{\II}_{\sigma} (\bxx, \xi,\lef(i),t^*_{n+1})\right],
\quad
\varphi^{\II}_{\sigma} (\bxx, \xi,i,t^*_{n+1}),
\quad
\eta\left[\varphi^{\II}_{\sigma} (\bxx, \xi,\rig(i),t^*_{n+1})\right],
$$ 
depend only on  $M^*_{n+1}$, on $c^*_{n+1}$, and on
$$   
\varrho\left[\varphi^{\II}_{\sigma} (\bxx, \xi,\lef(i),t^*_{n})\right],
\quad
\varphi^{\II}_{\sigma} (\bxx, \xi,i,t^*_{n}),
\quad
\eta\left[\varphi^{\II}_{\sigma} (\bxx, \xi,\rig(i),t^*_{n})\right].
$$
To this aim, we first assume that $c^*_{n+1}=i$ and we examine the value of the flag $f^*_{n+1}$. 
\begin{itemize}
\item
If $f^*_{n+1}=\typeu$, $\typev$ or $\typew$, by the very construction of the process, 
$\varphi^{\II}_{\sigma} (\bxx, \xi,i,t^*_{n+1})$ is determined
by $\varphi^{\II}_{\sigma} (\bxx, \xi,i,t^*_{n})$, $f^*_{n}$ and $z^*_{n}$. 
\item If $f^*_{n+1}=\typer$ or $\typeq$, 
$\varphi^{\II}_{\sigma} (\bxx, \xi,i,t^*_{n+1})$ is determined by  
$\varphi^{\II}_{\sigma} (\bxx, \xi,i,t^*_{n})$, by $z^*_{n}$, and by 
the knowledge of whether or not the sequence accepts the \typer\ or \typeq\ substitutions to $z^{*}_{n}$
at  site $i$. In turn, this only depends on 
$$\varrho\left[\varphi^{\II}_{\sigma} (\bxx,
\xi,\lef(i),t^*_{n})\right],
\quad
\varphi^{\II}_{\sigma} (\bxx, \xi,i,t^*_{n}),
\quad
\eta\left[\varphi^{\II}_{\sigma} (\bxx, \xi,\rig(i),t^*_{n})\right].
$$
\end{itemize}
This settles the case when $c^*_{n+1}=i$.
Now we assume that $c^*_{n+1}= \lef(i)$ and we examine the value of the flag $f^*_{n+1}$. 
\begin{itemize}
\item If $f^*_{n+1}=\typeu$, $\typev$ or $\typew$,  $\varphi^{\II}_{\sigma} (\bxx, \xi,\lef(i),t^*_{n+1})$ 
depends on $z^*_{n}$ and on whether the nucleotide
 $\varphi^{\II}_{\sigma} (\bxx, \xi,\lef(i),t^*_{n})$ is a purine or a pyrimidine.
This information is provided by
$\varrho\left[\varphi^{\II}_{\sigma} (\bxx, \xi,\lef(i),t^*_{n})\right]$.
\item  If $f^*_{n+1}=\typer$ or $\typeq$ and $z^*_{n}=C$ or $z^*_{n}=T$, 
 whether the sequence accepts the substitutions of types $\typer$ or $\typeq$
at site $\lef(i)$ depends only on  
$\varrho\left[\varphi^{\II}_{\sigma} (\bxx,
\xi,\lef(i),t^*_{n})\right]$ and 
  $\varphi^{\II}_{\sigma} (\bxx, \xi,i,t^*_{n})$.
\item If $f^*_{n+1}=\typer$ or $\typeq$ and $z^*_{n}=A$ or $z^*_{n}=G$,
 we observe that the corresponding substitution of types $\typer$ or $\typeq$ at site $\lef(i)$ could only turn an 
$A$ to a $G$ or vice versa. This does not affect the value of   
 $\varrho\left[\varphi^{\II}_{\sigma} (\bxx, \xi,\lef(i),t^*_{n+1})\right]$, which must be equal to 
$\varrho\left[\varphi^{\II}_{\sigma} (\bxx, \xi,\lef(i),t^*_{n})\right]$.
\end{itemize}
This settles the case when $c^*_{n+1}=\lef(i)$.
Since symmetric arguments hold when  $c^*_{n+1}= \rig(i)$, this concludes the proof.
\end{proof}

\begin{proof}[\ of proposition~\ref{p.r}]
Let $\prt^{a-1,b+1}(\bxx_{k}):=(x_{i})_{a-1\le i\le b+1}$. The result 
follows from proposition~\ref{p.i} since, for every $i$ in $\{a,  \ldots, b \}$,
$$\varphi^{c_{1},d_{1}}_{\sigma} (\bxx_{1}, \xi,i,s)\quad \mbox{and}\quad
\varphi^{c_{2},d_{2}}_{\sigma} (\bxx_{2},\xi,i,s)
$$ 
are both equal to
$$ \Theta_{\sigma,s}(x_{i-1}, x_{i}, x_{i+1} ,\xi_{i-1}[\sigma,s],  \xi_{i}[\sigma,s] ,
\xi_{i+1}[\sigma,s]).$$
\end{proof}

\section{Construction on the integer line}
\label{s.tpoaz}
The construction of the process on the integer line $\Z$ relies crucially on a
consequence of proposition~\ref{p.r} above, namely, the fact that, 
for every initial sequence $\bxx$, every site $i$ in $\Z$,
 every source of moves $\xi$, and  every
couple of times $s\ge\sigma$, 
the value of 
$$
\varphi^{a,b}_{\sigma}(\prt^{a,b}(\bxx),\xi,i,s)
$$
does not depend on $a$ and $b$ as soon as 
$$
a+1
\le i \le b-1.
$$ 
Hence the projective limit of the
system $(\varphi^{a,b}_{\sigma})_{a,b}$ when $a\to-\infty$ and
$b\to\infty$ exists trivially and defines a 
mesurable map 
$$
\Phi_{\sigma} \ : \ \A^{\Z} \times \Omega_{1}  \to  \D(\sigma, \A^{\Z}).
$$
Some previous observations about $\varphi^{a,b}_{\sigma}$ translate
immediately to $\Phi_{\sigma}$. 
\begin{prop}\label{p:redit}
{\bf (1)}
For every $s\ge t\ge\sigma$,
$$\Phi_{\sigma} (\bxx, \xi)(s) =  \Phi_{t}\left(   \Phi_{\sigma} (\bxx, \xi)(t) , \xi \right)(s). 
$$
{\bf(2)}
For every $s \ge  0$ and $t$,
$$
\Phi_{\sigma+t} (\bxx, \xi+t)(\sigma+t+s) = \Phi_{\sigma} (\bxx, \xi)(\sigma+s).$$
{\bf(3)}
Finally, $\Phi_{\sigma} (\bxx, \xi)(s)$ depends on $\xi$ only through $\xi[\sigma,s]=(\xi_{i}[\sigma,s])_{i \in \Z}$.
\end{prop}
For every $\bxx$ in $\A^{\Z}$, let $\P_{\bxx}$ denote the distribution of
$$
s \mapsto \Phi_{0}(\bxx, \cdot) (s),
$$ viewed as a random variable on $(\Omega_{1},\F_{1},\Q)$ with values
in $\D(0, \A^{\Z})$.  Then, as can be checked from
propositions~\ref{p.r} and \ref{p:redit} using the translation
invariance of Poisson processes, the family
$$
\left\{\,\P_{\bxx},\,\bxx\in\A^{\Z}\,\right\}
$$
defines a Feller Markov process in the sense of Liggett~\cite[chapter~1]{Lig}. From now on, we use the notation 
$$
X^{\bxx}_{s}=
\Phi_{0}(\bxx,
\cdot) (s), 
$$
and we sometimes omit the initial condition $\bxx$ of the Markov process $(X_{s})_{s \ge  0}$.

It is straightforward to check that the
construction in \cite{Lig}, based on infinitesimal generators, yields the same process. Let us call $\mathcal{G}$ the infinitesimal 
generator of the process yielded by the construction in \cite{Lig},
then $\mathcal G$ is well-defined and explicitly known at least for the 
Lipschitz functions on $\A^\Z$.

Let $(\widetilde{X}_{\bxx}^{\II}(s))_{s \geq 0}$ denote the Markov
process on $\A^{\Z}$ defined by 
$$\widetilde{X}_{\bxx}^{\II}(s)_{i} =
\left\{\begin{array}{ll}X_{\bxx}^{\II}(s)_{i}& \mbox{if}\ i \in \II,
\\
A & \mbox{else.}\end{array}\right.
$$  
By definition, $\widetilde{X}_{\bxx}^{\II}(s)$ converges to
$X^{\bxx}_{s}$ as $\II \to \Z$.  On the other hand, since the process
$(\widetilde{X}_{\bxx}^{\II}(s))_{s \geq 0}$ involves moves only on
the finite set of sites $\II$, its infinitesimal generator
$\mathcal{G}_{\II}$ can be readily computed. Moreover, 
$\mathcal{G}_{\II}$ converges to $\mathcal{G}$ as $\II \to \Z$,
at least on the set of real valued Lipschitz functions defined on
$\A^{\Z}$.  This is enough to identify $(X_{s})_{s \ge 0}$ with the
process yielded by the construction of \cite{Lig}, according to
Corollary 3.14 in \cite{Lig}.

A consequence of proposition~\ref{p.i} above is the following result.
\begin{prop}\label{p.i2}
For every sequence $\bxx:=(x_{i})_{i\in\Z}$, every
integer $i$ and every couple of times $s\ge\sigma$, 
$\Phi_{\sigma} (\bxx, \xi,i,s)$ depends
on $\bxx$
 and $\xi$ only through  
$$x_{i-1},\quad
x_{i},\quad
x_{i+1},\quad
\xi_{i-1}[\sigma,s],\quad
\xi_{i}[\sigma,s],\quad
\xi_{i+1}[\sigma,s].
$$ 
Indeed, using the function $\Theta_{\sigma,s}$ defined in
proposition~\ref{p.i},
$$
\Phi_{\sigma} (\bxx, \xi,i,s)    =    \Theta_{\sigma,s}(x_{i-1}, x_{i}, x_{i+1} ,\xi_{i-1}[\sigma,s], 
\xi_{i}[\sigma,s] ,
\xi_{i+1}[\sigma,s]).
$$
\end{prop}
A simple consequence of proposition~\ref{p.i2}
is the following proposition.
\begin{prop}\label{p.dat3}
Fix $\bxx$ and some subsets $I$ of the integer line at distance at least $3$ from each other, that is, such that
for every such pair $I\neq I'$ of subsets and every sites $i$ in $I$ and $i'$ in $I'$, $|i-i'|\ge3$.
Then,
the collections $\left[(X^{\bxx}_{s
\ge  0})_{i }\right]_{i
\in I}$ are independent.
\end{prop}
Here are some simple examples.
\begin{itemize}
\item
The collections $\left[(X_{s
\ge  0})_{i }\right]_{i
\le 0}$ and $\left[(X_{s
\ge0})_{i }\right]_{i
\ge3}$ are independent from each other.
\item
The collection 
$
\left[(X_{s
\ge  0})_{4i},(X_{s
\ge  0})_{4i+1}\right]_{i\in\Z}
$ 
is made of i.i.d.\ random variables with values in $\A^2$.
\item
The remaining values form a collection 
$
\left[(X_{s
\ge  0})_{4i+2},(X_{s
\ge  0})_{4i+3}\right]_{i\in\Z},
$ 
which is also made of i.i.d.\ random variables with values in $\A^2$,
with the same distribution as the collection in the preceding item,
while being not independent from it.
\item
The collection
$
\left[(X_{s
\ge  0})_{3i}\right]_{i\in\Z},
$
is made of i.i.d.\ random variables with values in $\A$.
\item
For every $\ell\ge0$ and every $\delta\ge\ell+3$,
the collection
$
\left[(X_{s
\ge  0})_{\{\delta i,\ldots,\delta i+\ell\}}\right]_{i\in\Z},
$
is made of i.i.d.\ random variables with values in $\A^{\ell+1}$.
\end{itemize}
\begin{prop}\label{p.st}
There exists a unique stationary distribution of $(X_{s})_{s\ge0}$ on $\A^{\Z}$. 
\end{prop}
\begin{definition}\label{d.mu}
Let $\mu$ denote the stationary distribution  of $(X_{s})_{s\ge0}$.
Let $\mu_{a,b}$ denote the measure $\mu_\II$ whose existence is ensured by
proposition~\ref{p.ef},
when $\II=\ab$.
\end{definition}
A consequence of the uniqueness of $\mu$ is its invariance by the
transformations that preserve the dynamics, for instance the
translations of $\Z$. Likewise, $\mu_{a,b}$ is invariant with respect to the translations of the discrete circle $ \ab \cong \Z / (b-a+1) \Z$.
Moreover, if $\ab$ is the image of
$\{c,\ldots,d\}$ by a translation of $\Z$, $\mu_{a,b}$ coincide with
 the image of $\mu_{c,d}$
by this translation. Other properties are in proposition~\ref{p.st2}.
\begin{prop}\label{p.st2}
{\bf(1)}
Assume that the integers  $a$, $b$, $c$ and $d$ are such that  
$$
c+1 \le a\le b \le d-1.
$$
Then,
$$
\prt^{a,b}(\mu_{c,d}) = \prt^{a,b}(\mu).
$$ 
{\bf(2)}
The Markov process $(X_{s})_{s\ge0}$ is ergodic. That is, for every initial condition
$\bxx$ in $\A^{\Z}$, $X^{\bxx}_{s}$ converges in distribution to $\mu$ as $s\to\infty$.  
\\
{\bf(3)}
Finally, the independence properties stated above hold with respect to $\mu$ as well.
\end{prop}
The considerations above justify the following definition.
\begin{definition}[Stationary frequencies]\label{d.ff}
Define the stationary frequency $F(x_1\cdots x_k)$ of every polynucleotide $x_1\cdots
x_k$ as
$$
F(x_1\cdots x_k):=\mu_{0,k+1}(\A\times\{(x_1,\ldots,x_k)\}\times\A).
$$
\end{definition}
Hence, $F(x_1\cdots x_k)$ is also
$$
F(x_1\cdots x_k)
=\mu\left(\A^{\{\ldots,-1,0\}}\times\big\{(x_1,\ldots,x_k)\big\}\times\A^{\{k+1,k+2,\ldots\}}\right).
$$

Due to the independence properties of $\mu$ stated in proposition~\ref{p.st2}, $\mu$ is clearly ergodic with respect to the translations in $\Z$, so we may
as well define stationary frequencies as:
$$F(x_1\cdots x_k) = \lim_{a \to -\infty, b \to +\infty}  \frac{1}{b-a+1}\sum_{i=a}^{b} \un(X_{i+1} \cdots X_{i+k} =  x_1\cdots x_k)$$
where the distribution of $(X_{i})_{i \in \Z}$ is $\mu$, $\un(A)$ denotes the indicator function of the event  $A$, and the above limit holds in the almost sure sense.

\begin{proof}[\ of propositions~\ref{p.st} and \ref{p.st2}]
Fix $a$ and $b$ such that $a \le b$.  From proposition~\ref{p.r}, for
every $\bxx$ in $\A^{\Z}$ and every $s \ge 0$,
$$\prt^{a,b}\left(X^{\bxx}_{s}\right) =
\prt^{a,b}\left(X^{a-1,b+1}_{\bxx}(s)\right).$$ Now, according to
proposition~\ref{p.ef}, $X^{a-1,b+1}_{\bxx}(s)$ converges in
distribution towards $\mu_{a-1,b+1}$ as $s\to\infty$. As a
consequence, $\prt^{a,b}\left(X^{\bxx}_{s}\right)$ converges to
$\prt^{a,b}(\mu_{a-1,b+1})$. Proposition~\ref{p.r} again shows that
$$
\left( \prt^{a,b}(\mu_{a-1,b+1})  \right)_{a\le b} 
$$ 
is a coherent family of
probability distributions. 
Hence there exists a unique probability distribution  $\mu$ on $\A^{\Z}$ such that
$\prt^{a,b}(\mu) =\prt^{a,b}(\mu_{a-1,b+1})$.
  The convergence of  $\prt^{a,b}\left(X^{\bxx}_{s}\right)$ in distribution towards  $\prt^{a,b}(\mu)$  
for every $a\le b$ implies
that, for every $\bxx$ in $\A^{\Z}$, 
$X^{\bxx}_{s}$ converges in distribution to $\mu$. 
The existence and the uniqueness of an invariant distribution, equal to
$\mu$, follow easily. The other properties are obvious.
\end{proof}

\section{R/Y encodings}
\label{s.tryp}
The situation of the R/Y process is even simpler.
\subsection{On finite intervals}
\begin{definition}
\label{d.ryq}
Introduce the R/Y configuration at time $s$ as the $\{R,Y\}$-valued
function
$$Z^{\II}_{\bxx}(s):=\pi(\varphi_{\sigma}^{\II}(\bxx,\xi,i,s)).
$$
\end{definition}
\begin{prop}
\label{p.ry}
For every sequence $\bxx:=(x_{i})_{i\in\II}$, every
site $i$ in $\II$ and every time $s\ge \sigma$,
the value of $Z^{i}_{\bxx}(s)$, which depends a
priori on all the information in $\bxx$ and $\xi$, is in fact 
measurable with respect to
$\prp(x_i)$ and $\xi[\sigma,s]$.
\end{prop}
\begin{proof}[\ of proposition~\ref{p.ry}]
The substitutions associated to clocks of types \typew, \typer\ and \typeq\ do not change
the value of $Z^{\II}_{\bxx}(s)$. Hence the moves of the function 
$s\mapsto Z^{\II}_{\bxx}(s)$
are determined by the clocks $\UU_i$ and $\VV_i$. Erasing the $\RR_i$ and $\QQ_{i}$
clocks is like setting every $r^{y}_z$ to $0$, in which case the sites
evolve independently. This proves the proposition.
\end{proof}
\begin{corollary}
\label{c.ry}
For every collection of sites $\II$, the coordinates of
$$
Z^{\II}_{\bxx}(s)$$ 
are independent and
each converge in distribution, when $s\to\infty$, to the distribution
on $\{R,Y\}$ equal to
$t_Y\,\delta_Y+t_R\,\delta_R$, that is,
$$
\P^{\II}_{\bxx}((Z_s)_i=Y)\to t_Y,\quad
\P^{\II}_{\bxx}((Z_s)_i=R)\to t_R.
$$
\end{corollary}
We recall from theorem~\ref{t.ry} in section~\ref{s.idotr} that
$t_Y+t_R=1$ and
$$
t_Y:=\frac{v_C+v_T}{v_A+v_T+v_C+v_G},\quad
t_R:=\frac{v_A+v_G}{v_A+v_T+v_C+v_G}.
$$

\subsection{On the full line}
Proposition~\ref{p.ry2} below is a straightforward consequence and 
equivalent to theorem~\ref{t.ry} in section~\ref{s.idotr}.
\begin{prop}
\label{p.ry2}
There exists a unique stationary distribution of $(Z_s)_{s\ge0}$ on
$\{R,Y\}^{\Z}$. This measure is the product mesure $\nu^{\otimes\Z}$, where
$$
\nu(Y):=t_Y,\quad
\nu(R):=t_R.
$$
\end{prop}
Theorem~\ref{t.ry} yields relations, which hold irrespective of the values of the mutation rates $r^y_{x}$,
namely,
\begin{eqnarray*}
F(CG)+F(CA)+F(TG)+F(TA) &  = & t_{Y}\,t_{R},\\
F(CC)+F(CT)+F(TC)+F(TT) &  = &  t_{Y}^{2},\\
F(AC)+F(AT)+F(GC)+F(GT) &  = &  t_{R}\,t_{Y},\\
F(CC)+F(CT)+F(TC)+F(TT) & = &  t_{R}^{2}.
\end{eqnarray*}
These relations are simple enough to write.  However, they yield
awkward formulas for the individual frequencies of nucleotides and
dinucleotides, in full generality.

We now comment on theorem~\ref{t.ry} and proposition~\ref{p.ry2}.  
This phenomenon is a
consequence of the construction with exponential clocks, where the
locations of the YpR dinucleotides coincide, before and after the
superposition of the exponential clocks which represent the YpR
mutations.  Hence the overall frequency of these $4$ dinucleotides
should be the same for every values of the rates $r^{x}_{y}$. Indeed
the right hand side is the value of their overall frequency when
$r^x_{y}=0$ for every $xy$. Also, this can be checked directly from
the equations which describe the equilibrium of these $4$
dinucleotides. Finally, the product $t_{R}\,t_{Y}$ in the right hand
side coincides with the product of the frequencies of the pyrimidines
and of the purines, that is, the product
$$
F(Y)\,F(R)=(F(C)+F(T))\,(F(A)+F(G)).
$$ Once again, this can be seen on the construction with exponential
clocks.  Property (a) means that one can fuse $C$ and $T$ into a
single state $Y$ (pyrimidine), and $G$ and $A$ into a single state $R$
(purine).  The YpR substitutions have no effect on $R$ and $Y$,
hence the sites are i.i.d.  Each pyrimidine mutates at rate 
$$s_{R}:=v_A+w_T+w_C+v_G,
$$
to a purine with probability $(v_{A}+v_{G})/s_{R}$ and to a pyrimidine
otherwise.  Each purine mutates at rate 
$$s_{Y}:=w_A+v_T+v_C+w_G,
$$
to a pyrimidine with
probability $(v_{T}+v_{C})/s_{Y}$ and to a purine otherwise. The net
effect is that every pyrimidine becomes a purine at rate $v_{A}+v_{G}$
and that every purine becomes a pyrimidine at rate $v_{C}+v_{T}$,
hence the stationary measure $F(Y):=F(C)+F(T)$ of the pyrimidines is
proportional to $v_{C}+v_{T}$ and the stationary measure
$F(R):=F(A)+F(G)$ of the purines is proportional to $v_{A}+v_{G}$.

\section*{Part B Computation}
\addcontentsline{toc}{section}{Part B Computation}

\section{General case}
\label{s.tgc}

\subsection{Polynucleotide frequencies}
\label{ss.pf}
From the construction given in part A, knowing the stationary
distribution of the Markov process $(X^{\II}(s))_{s \geq 0}$ with
$\II= \{ a-1,\ldots, b+1 \}$ is enough to compute
$\Pi^{a,b}(\mu)$. Since $(X^{\II}(s))_{s \geq 0}$ lives on the state
space $\A^{\II}$, computing its stationary distribution amounts to
solving a linear system of size $\# \A^{\II} \times \# \A^{\II}$. 
Computing the equilibrium frequency of polynucleotides of length $N$ thus
requires solving a $4^{N+2} \times 4^{N+2}$ linear system.

Theorem~\ref{t.ration}
in section~\ref{s.idotr} follows from these considerations and from
Cram\'er's formula.
However, even moderate lengths of polynucleotides
($N=4$, say) lead to fairly large
linear systems, so finding ways of lowering the dimension of the
system to be solved is a critical issue, if solvability of the model
is to be considered something more than a mere theoretical
possibility.

For single nucleotides and the restricted class of YpR dinucleotides,
an autonomous linear subsystem can be isolated, and this is discussed
in section~\ref{ss.computYpR} below.

For general polynucleotides, symmetries can be used to reduce the
computational burden. Using the invariance of $(X^{\II}(s))_{s \geq
0}$ with respect to translations on the discrete circle, the linear
system yielding the stationary distribution of $(X^{\II}(s))_{s \geq
0}$ can be reduced to a linear system of size
$$m(N+2) \times m(N+2),
$$
where $m(k)$ denotes the number of distinct orbits in
$\A^{k}$ under translations. It is a well-known counting
result, see e.g.~\cite{Rot}, that   
$$m(k)= \frac{1}{k} \sum_{d | k} \phi(d)\,4^{k/d},
$$ where $\phi$ stands for Euler's indicator
(the number of primes to a number). Hence, when $k\to\infty$,
$$
m(k)\sim4^k/k.
$$
This 
remark also achieves significant improvement for small values of $N$,
see the table below.
The last columns give the value of $1/(N+2)$ and the exact ratio 
 $m(N+2)/4^{N+2}$ of the sizes of the
two linear systems.
\begin{center}
\begin{tabular}{|c|c|c|c|c|}
\hline
$N$ & $4^{N+2}$ & $m(N+2)$ & $1/(N+2)$ & Ratio
\\
\hline
$2$ & $256$ & $70$ & $25\%$ & $27.3 \%$
\\
$3$ & $1024$ & $208$ & $20\%$ & $20.3 \%$
\\
$4$ & $4096$ & $700$ & $16.7\%$ & $17.1 \%$
\\
\hline
\end{tabular}\end{center}
\begin{remark}
One could think of the following alternative reduction
to lower the size of the linear systems.
From the results of part~A, to find $\prt^{a,b}\left( X^{\II}(s) \right)$
with $\II= \{ a-1,\ldots, b+1 \}$, it is sufficient to find the
invariant 
distribution of the Markov chain
$$\left( \prr(X^{\II}(s)_{a-1}) , \prt^{a,b} \left( X^{\II}(s) \right)
,\pry(X^{\II}(s)_{b+1}) \right).
$$ 
This chain lives on a state space of size $4^{N} \times 3^{2}$.  
Hence, this remark reduces the size
of the system by a factor $9/16\approx56 \%$. On the other hand, 
the translation invariance is
lost.
All in all, using the translation
invariance described above yields more effective reductions.
\end{remark}

Two approches to computing equilibrium frequencies of polynucleotides
may be considered. One can solve the linear system numerically with fixed
values of the parameters (with finite or infinite precision
arithmetic), within reasonable time for $N\le4$.
One can also solve this symbolically, 
a task that we performed only for $N=2$, using
a restricted version of the model possessing a single free parameter,
and additional symmetries, see section~\ref{s.soa}.

\subsection{Nucleotide and YpR
  dinucleotide frequencies}
\label{ss.computYpR}

Recall that $\fr(x_{1}\cdots x_{k})$, introduced formally in
 definition~\ref{d.ff} in section~\ref{s.tpoaz},
 denotes the stationary frequency of the polynucleotide $x_{1}\cdots x_{k}$. 
In this section, we show how to compute 
$\fr(x)$ for every nucleotide $x$,
and $\fr(xy)$ for every YpR dinucleotide $xy$.
\begin{definition}\label{d.freq}
Introduce
$$
\fr(Y):=\fr(C)+\fr(T),\quad\fr(R):=\fr(G)+\fr(A).
$$
Similar conventions are valid for polynucleotides, for instance
$$
\fr(YR):=\sum_{\prp(x)=Y}\sum_{\prp(y)=R}\fr(xy).
$$
Likewise,
$$
\fr(Yy):=\fr(Cy)+\fr(Ty),\quad
\fr(xR):=\fr(xA)+\fr(xG).
$$
\end{definition}
We introduce some notations, related to the rates of the simple substitutions.
\begin{definition} \label{d.tyr}
For every nucleotide $x$, let $s_{x}$ denote the sum of the rates of mutations from $x$, hence
$s_{A}:=s_{R}=:s_{G}$ and $s_{C}:=s_{Y}=:s_{T}$, with
$$
s_{R}:=w_{A}+v_{T}+v_{C}+w_{G},
\quad
s_{Y}:=v_{A}+w_{T}+w_{C}+v_{G}.
$$
Likewise, let 
$$
u_x:=v_x-w_x,\quad
v:=\sum_{x}v_x,\quad
w:=\sum_{x}w_{x}.
$$
Finally, let $t_{A}:=t_R=:t_G$ and $t_C:=t_Y=:t_T$, with $t_{R}+t_{Y}=1$ and
$$
t_{R}:=(v_{A}+v_{G})/v,
\quad
t_{Y}:=(v_{T}+v_{C})/v.$$
\end{definition}
We turn to some notations related to the effect of the $r_{x}^{y}$ 
substitutions on nucleotides.
\begin{definition}\label{d.pxy}
For every nucleotide $x$ and every YpR dinucleotide $yz$, 
introduce $p_{yz}(x)$ as 
the rate at which $x$ appears (or disappears if $p_{yz}(x)$ is
negative) because of
 the dinucleotide substitutions 
associated to $yz$. Hence,
$$
\begin{array}{ccc}
p_{CG}(T):=r^G_{T}=:-p_{CG}(C),
& \quad &
p_{CG}(A):=r^C_{A}=:-p_{CG}(G),
\\
p_{TA}(C):=r^{A}_{C}=:-p_{TA}(T),
& \quad &
p_{TA}(G):=r^{T}_{G}=:-p_{TA}(A),
\\
p_{CA}(T):=r^A_{T}=:-p_{CA}(C),
& \quad &
p_{CA}(G):=r^C_{G}=:-p_{CA}(A),
\\
p_{TG}(A):=r^{T}_{A}=:-p_{TG}(G),
& \quad &
p_{TG}(C):=r^{G}_{C}=:-p_{TG}(T).
\end{array}
$$
Finally, for every $x$ and every $yz$ which is not a YpR dinucleotide, let
$$
p_{yz}(x):=0.
$$
\end{definition}
Equilibrium for the nucleotides yields the following relations.
\begin{prop}\label{p.un}
For every nucleotide $x$,
$$
s_{x}\,\fr(x)=v_{x}-u_x\,t_{x}+\sum_{yz} p_{yz}(x)\,\fr(yz).
$$
Furthermore, $\fr(R)=t_{R}$ and $\fr(Y)=t_{Y}$.
\end{prop}
Note that the values of $\fr(R)$ and $\fr(Y)$ are independent of the YpR substitution
rates $r^y_{x}$. The underlying reason for
this a priori surprising fact is in proposition~\ref{p.ry2}.
The proof of proposition~\ref{p.un} is in section~\ref{ss.pro}.

From proposition~\ref{p.un}, the values of $F(xy)$  for
every YpR  dinucleotide $xy$ determine $F(z)$ for every nucleotide $z$. 
To compute $F(xy)$ for these $4$ dinucleotides $xy$,
we need some notations related to the effects of the $r_{x}^{y}$ mutations on dinucleotides.
\begin{definition}\label{d.pzt}
For every couple of YpR dinucleotides $xy$ and $zt$, 
introduce $p_{zt}(xy)$ as 
the rate at which $xy$ appears (or disappears if $p_{zt}(xy)$ is
negative), 
due to the existence of the
dinucleotide $zt$. Hence, assuming that $\{x,x^\ast\}=\{C,T\}$ and 
that $\{y,y^\ast\}=\{A,G\}$,
$$
p_{xy}(xy):=-r^{x}_{y^\ast}-r^{y}_{x^\ast},
\quad
p_{xy}(xy^\ast):=r^{x}_{y^\ast},
\quad
p_{xy}(x^\ast y):=r^{y}_{x^\ast},
\quad
p_{xy}(x^\ast y^\ast):=0.
$$
For instance,
$$
p_{CG}(CG):=-r^G_{T}-r^C_{A},
\quad
p_{CA}(CG):=r^C_{G},
\quad
p_{TG}(CG):=r^G_{C},
\quad
p_{TA}(CG):=0.
$$
Finally, let
$$
q_{xy}:=-p_{xy}(xy)=r^{x}_{y'}+r^{y}_{x'}.
$$
\end{definition}
Equilibrium for the YpR dinucleotides yields the following relations. 
Recall the notations in definition~\ref{d.freq}.
\begin{prop}\label{p.deux}
For every YpR dinucleotide $xy$, 
$$
(v+w)\,\fr(xy) +u_x\,\fr(Yy)
+u_y\,\fr(xR) =
v_{x}\,\fr(y)+v_{y}\,\fr(x)+\sum_{zt}p_{zt}(xy)\,\fr(zt).
$$
\end{prop}
The proof of proposition~\ref{p.deux} is in section~\ref{ss.pro}.

From propositions~\ref{p.un} and \ref{p.deux},
the $8$ unknown frequencies we are looking for, namely the $4$ frequencies of the nucleotides 
and the $4$ frequencies of the YpR dinucleotides, solve a system
of $8$ linear equations.
One can show that the determinant of this system is not zero, hence the $8$ frequencies 
are entirely determined by this system.


A simpler way to proceed is to write the frequency of each nucleotide as 
an affine function of the $4$ fequencies of YpR dinucleotides, then to plug these expressions 
 in the $4$ last equations. We state this as theorem~\ref{t.slin}
 below, which precises theorem~\ref{t.sn} in section~\ref{s.idotr}.
\begin{definition}\label{d.vecteur}
Let $\fmat$ denote the $4\times1$ vector of the frequencies of the YpR dinucleotides, that is,
$$
\fmat:=\left(\begin{array}{c}F(CG)\\
F(CA)\\
F(TG)\\
F(TA)\end{array}\right).
$$
\end{definition}
\begin{theorem}[YpR frequencies]\label{t.slin}
The YpR frequencies solve a linear system
$$
((v+w)\,\mathrm{Id}+\umat+\rmat)\cdot\fmat=\vmat,
$$
where the $4\times 4$ matrices $\umat$ and $\rmat$ and the $4\times1$ vector $\vmat$ are defined 
below and depend on the substitution rates $v_{x}$, $w_{x}$ and $r^{y}_{x}$.
\end{theorem}
\begin{definition}\label{d.normfreq}
For every YpR dinucleotide, let
$$v^*_x:=v_{x}/s_{R},\qquad
v^*_{y}:=v_{y}/s_{Y}.$$
\end{definition}
The matrix $\umat$ is
$$
\umat:=\left(\begin{array}{cccc}u_{C}+u_{G}&u_{G}&u_{C}&0\\
u_{A}&u_{C}+u_{A}&0&u_{C}\\
u_{T}&0&u_{T}+u_{G}&u_{G}\\
0&u_{T}&u_{A}&u_{T}+u_{A}\end{array}\right).
$$
The coefficients of the matrix $\rmat$ are
$$
\rmat_{xy,zt}:=-p_{zt}(xy)-v^*_{x}\,p_{zt}(y)-v^*_{y}\,p_{zt}(x).
$$
Finally, the coefficients of the matrix $\vmat$ are
$$
\vmat_{xy}:=v_{x}\,\frac{v_{y}-u_{y}\,t_{y}}{s_{y}}+v_{y}\,\frac{v_{x}-u_{x}\,t_{x}}{s_{x}}.
$$
Note that every $\vmat_{xy}$ is positive.

We now write the coefficients of $\rmat$ more explicitly.
Each column $\rmat_{\cdot,xy}$ of $\rmat$ depends on the YpR rates of substitution through
$r^{x}_{y^\ast}$ and $r^{y}_{x^\ast}$ only, and through affine functions.
More precisely,
$$
\rmat_{xy,xy}=(1+v^{*}_{x})\,r^{x}_{y^\ast}+(1+v^{*}_{y})\,r^{y}_{x^\ast},
$$
Furthermore,
$$
\rmat_{xy^\ast,xy}=-(1+v^{*}_{x})\,r^{x}_{y^\ast}+v^{*}_{y^\ast}\,r^{y}_{x^\ast},
\quad
\rmat_{x^\ast y,xy}=-(1+v^{*}_{y})\,r^{y}_{x^\ast}+v^{*}_{x^\ast}\,r^{x}_{y^\ast},
$$
and
$$
\rmat_{x^\ast y^\ast,xy}=-v^*_{x^\ast}\,r^{x}_{y^\ast}-v^*_{y^\ast}\,r^{y}_{x^\ast}.
$$
For instance, the $CG$ column of $\rmat$ is
$$
\rmat_{\cdot,CG}:=\begin{array}{cc}
\begin{array}{c}CG\\ CA\\ TG\\ TA\end{array}
 &
 \!\!\!\!
\left(\begin{array}{c}r^C_{A}+r^G_{T}+v_{C}r^C_{A}/s_{R}+v_{G}r^G_{T}/s_{Y}\\
-r^C_{A}-v_{C}r^C_{A}/s_{R}+v_{A}r^G_{T}/s_{Y}\\
-r^G_{T}-v_{G}r^G_{T}/s_{Y}+v_{T}r^C_{A}/s_{R}\\
-v_{T}r^C_{A}/s_{R}-v_{A}r^G_{T}/s_{Y}\end{array}\right).
\end{array}
$$

\subsection{Proofs}
\label{ss.pro}
\begin{proof}[\ of proposition~\ref{p.un}]
Assume that the distribution of $(X_{i})_{i \in \Z}$ is the stationary
measure $\mu$ introduced in definition~\ref{d.mu} in section~\ref{s.tpoaz}.
By the definition of the dynamics, equilibrium for nucleotide $x$ at
site 
$i$ reads
\begin{eqnarray*} 
s_{x}\,\P(X_{i}=x)
& = &    
\sum_{\prp(z)=\prp(x)} w_{x}\,\P(X_{i} = z )  +  
\sum_{\prp(z)\neq\prp(x)} v_{x}\,\P(X_{i} = z )
\\
& &
\hspace{5ex}+\sum_{\prp(y)=\prp(x)\neq\prp(z)} p_{yz}(x)\,\P(X_{i}X_{i+1} = yz) 
\\
& &
\hspace{5ex}+\sum_{\prp(y)\neq\prp(x)=\prp(z)} p_{yz}(x)\,\P(X_{i-1}X_{i} = yz).
\end{eqnarray*}
Note that one of the last two sums in the expression above is always zero, which one depending on 
whether $x$ is a purine or a pyrimidine. 

Using the translation invariance of the stationary distribution, and extracting constant 
factors from sums, this reduces to
\begin{equation}\label{e.reduite} 
s_{x}\,\fr(x)\,\fr(yz)=
\sum_{yz} p_{yz}(x) +
w_{x} \sum_{\pi(z)=\pi(x)}   \fr(z)   + v_{x} \sum_{\pi(z)\neq\pi(x)}  \fr(z).
\end{equation}
Using the fact that $p_{yz}(x) = -p_{yz}(x^\ast)$, and summing, on the one hand the above 
equilibrium equations for $x=A$ and $x=G$, and on the other hand for $x=C$ and $x=T$, 
we obtain a linear system of two equations involving $F(R)$ and
$F(Y)=1-F(R)$.
For instance,
$$
s_{R}\,\fr(R) = (w_{A}+w_{G})\,\fr(R) + (v_{A}+v_{G})\,\fr(Y).
$$
Solving this for $\fr(R)$ and $\fr(Y)$ yields the first assertion of the proposition, 
plugging these values into equation~(\ref{e.reduite}) yields the second assertion.
\end{proof}
\begin{proof}[\ of proposition~\ref{p.deux}]
As in the proof of proposition~\ref{p.un} above, we assume that the distribution of $(X_{i})_{i \in \Z}$ is 
the stationary measure $\mu$.
Equilibrium for a YpR dinucleotide $xy$ located
at the pair of sites $(i, i+1)$ reads as the equality of the exit rate
and the entrance rate. Both are due to single substitutions and to
double substitutions. The exit rate of single substitutions has size
$$
(s_{x}+s_{y})\,\P(X(i)X(i+1) = xy).
$$
The entrance rate due to the single transitions has size
$$
\sum_{\prp(z)=\prp(x)} w_{x} \,  \P(X(i)X(i+1) = zy)   +
\sum_{\prp(t)=\prp(y)} w_{y} \,  \P(X(i)X(i+1) = xt).
$$
The entrance rate due to the single transversions has size
$$
 \sum_{\prp(z)\neq\prp(x)} v_{x}\,   \P(X(i)X(i+1) = zy) + 
 \sum_{\prp(t)\neq\prp(y)} v_{y}\,   \P(X(i)X(i+1) = xt). 
$$
Finally, the rate of double substitutions, counted as an entrance
rate, has  size
$$
    \sum_{\prp(z)=\prp(x)} p_{zy}(x) \,\P(X(i) X(i+1) = zy)    +
\sum_{\prp(t)=\prp(y)} p_{xt}(y)\, \P(X(i) X(i+1) = xt).
$$
An important point to notice is that no YpR mutation affecting the pairs of sites $(i-1,i)$ or $(i+1,i+2)$ can occur when $X(i) X(i+1) = xy$ or lead to $X(i) X(i+1)=xy$, since $x$ is a pyrimidine and $y$ is a purine.  
This fact rules out probabilities of trinucleotides from appearing in the above equation, which 
we could not avoid if $xy$ was not an YpR dinucleotide.

Using the facts that 
$$
\sum_{\prp(z)\neq\prp(x)} \P(X(i)X(i+1) = zy) = 
\P(X(i+1)=y) -   \sum_{\prp(z)=\prp(x)}   \P(X(i)X(i+1) = zy),
$$ 
and that
$$\sum_{\prp(t)\neq\prp(y)}   \P(X(i)X(i+1) = xt) = 
\P(X(i)=x) -\sum_{\prp(t)=\prp(y)}    \P(X(i)X(i+1) = xt),$$ 
and the translation invariance of the stationary distribution, 
we obtain the identity stated in proposition~\ref{p.deux}.
\end{proof}

\section{Uniform simple rates}
\label{s.ussr}
In this section, we study the influence of the rates $r^{y}_{x}$ of YpR substitution rates,
looking at the case 
when,  for every nucleotide $x$, 
$$
v_{x}=w_{x}=1.
$$ 
Then, for every YpR dinucleotide $xy$,
$$
v=w=s_{x}=s_{y}=4,\quad u_{x}=0,\quad \vmat_{xy}=1/2,\quad v^*_{x}=v^*_{y}=1/4.
$$
We assume, as a further simplification, that the YpR substitution rates from CpA and from
TpG are zero, that is,
$$
r^C_{G}=r^{A}_{T}=r^G_{C}=r^T_{A}=0.
$$
Hence, $\rmat_{\cdot,CA}=\rmat_{\cdot,TG}=0$, and
$$
\rmat_{\cdot,CG}=\frac14\left(\begin{array}{c}5r^{C}_{A}+5r^{G}_{T}\\
-5r^{C}_{A}+r^{G}_{T}\\
r^{C}_{A}-5r^{G}_{T}\\
-r^{C}_{A}-r^{G}_{T}\end{array}\right),
\quad
\rmat_{\cdot,TA}=\frac14\left(\begin{array}{c}-r^{A}_{C}-r^{T}_{G}\\
-5r^{A}_{C}+r^{T}_{G}\\
r^{A}_{C}-5r^{T}_{G}\\
5r^{A}_{C}+5r^{T}_{G}\end{array}\right).
$$ 
Recall that $q_{CG}=r^{G}_{T}+r^{C}_{A}$ and $q_{TA}=r^{A}_{C}+r^{T}_{G}$.
\begin{definition}\label{d.dk}
Let 
$$
K:=32+5\,(q_{CG}+q_{TA})+3\,q_{CG}\,q_{TA}/4.
$$
\end{definition}
Our following result gives the frequencies of the
 YpR dinucleotides.
\begin{prop}\label{p.fypr}
Assume that the simple substitution rates are $1$ and that no YpR substitution occur from CpA or TpG.
Then, for every YpR dinucleotide $xy$,
$$
F(xy)=\frac1{16}\,\left(1+\frac{k(xy)}{K}\right),
$$
with 
\begin{eqnarray*}
k(CG) & := & q_{TA}-5\,q_{CG}-3\,q_{CG}\,q_{TA}/4,
\\
k(CA) & := & (5r^{C}_{A}-r^{G}_{T})\,(1+3\,q_{TA}/16)+(5r^{A}_{C}-r^{T}_{G})\,(1+3\,q_{CG}/16),
\\
k(TG) & := & (5r^{G}_{T}-r^{C}_{A})\,(1+3\,q_{TA}/16)+(5r^{T}_{G}-r^{A}_{C})\,(1+3\,q_{CG}/16),
\\
k(TA) & := & q_{CG}-5\,q_{TA}-3\,q_{CG}\,q_{TA}/4.
\end{eqnarray*}
\end{prop}
For instance,
$$
F(CG)=\frac{32+6\,q_{TA}}{16\,K},\qquad
F(TA)=\frac{32+6\,q_{CG}}{16\,K}.
$$ 
One can check that
$$
k(CG)+k(CA)+k(TG)+k(TA)=0.
$$
The frequencies of the nucleotides follow.
\begin{corollary}\label{c.fypr}
In the setting of proposition~\ref{p.fypr}, 
\begin{eqnarray*}
4\,F(A) & = & 1+r^{C}_{A}\,\frac{32+6\,q_{TA}}{16\,K}-r^{T}_{G}\,\frac{32+6\,q_{CG}}{16\,K},
\\
4\,F(G) & = & 1-r^{C}_{A}\,\frac{32+6\,q_{TA}}{16\,K}+r^{T}_{G}\,\frac{32+6\,q_{CG}}{16\,K},
\\
4\,F(C) & = & 1-r^{G}_{T}\,\frac{32+6\,q_{TA}}{16\,K}+r^{T}_{C}\,\frac{32+6\,q_{CG}}{16\,K},
\\
4\,F(T) & = & 1+r^{G}_{T}\,\frac{32+6\,q_{TA}}{16\,K}-r^{T}_{C}\,\frac{32+6\,q_{CG}}{16\,K}.
\end{eqnarray*}
\end{corollary}
If one assumes furthermore that $r^{T}_{C}=r^{T}_{G}=0$, that is, that CpG is the only active dinucleotide,
further simplifications occur.
\begin{prop}\label{p.fypr2}
Assume that the simple substitution rates are $1$ and that no YpR
substitution occur from CpA, TpG or TpA.
Then, for every YpR dinucleotide $xy$,
$$
F(xy)=\frac1{16}\,\left(1+\frac{k(xy)}{32+5(r^{G}_{T}+r^{C}_{A})}\right),
$$
with
$$
k(CG):=-5(r^{C}_{A}+r^{G}_{T}),\quad
k(CA):=5r^{C}_{A}-r^{G}_{T},
$$
and
$$
k(TG):=5r^{G}_{T}-r^{C}_{A},
\quad
k(TA):=r^{G}_{T}+r^{C}_{A}.
$$
\end{prop}
For instance,
$$
F(CA)=\frac1{16}\,\frac{32+4r^{G}_{T}+10r^{C}_{A}}{32+5r^{C}_{A}+5r^{G}_{T}}.
$$
\begin{corollary}\label{c.fypr2}
In the setting of proposition~\ref{p.fypr2}, 
$$
F(x)=\frac1{4}\,\left(1+\frac{k(x)}{32+5(r^{C}_{A}+r^{G}_{T})}\right),
$$
with
$$
k(T):=2r^{G}_{T},\quad
k(C):=-2r^{G}_{T},\quad
k(A):=2r^{C}_{A},\quad
k(G):=-2r^{C}_{A}.
$$
\end{corollary}



\section{Symmetric rates}
\label{s.wms}
\subsection{Models}
Consider  property (b) below.
\begin{quote}
\textbf{Property (b)}
The substitution rates respect the complementarity of the nucleotides.
\end{quote}
This means, first, that the rate of substitution from $x$ to $y$ and from $x^\ast$ to $y^\ast$ coincide, 
for every nucleotides $x$
and $y$, where we recall that the involution $z\mapsto z^\ast$ is 
defined by 
$$A^\ast:=T,\quad
T^\ast:=A,\quad
C^\ast:=G,\quad
G^\ast:=C.
$$
This means also that the rates of YpR substitutions
 from $CG$ to $CA$ and to $TG$ coincide, and that the rates of YpR substitutions
from $TA$ to $CA$ and to $TG$ coincide, that is,
$$
r^{C}_{A}=r^{G}_{T}=:r_{W},\quad r^{A}_{C}=r^{T}_{G}=:r_{S}.
$$
This means finally that the rates of YpR substitutions
 from $CA$ and from $TG$ to $CG$ coincide, and that the rates of YpR substitutions
from $CA$ and from $TG$ to $TA$ coincide

As regards the single substitutions, the most general model such that (a) and (b) hold is described by matrices
$$
\left(
\begin{array}{cccc}
\cdot & v_{W} & v_S & w_S
\\
v_{W} & \cdot & w_S & v_S
\\
v_{W} & w_W & \cdot & v_S
\\
w_W & v_{W} & v_S & \cdot
\end{array}
\right),
$$
where $v_S$, $v_W$, $w_S$ and $w_W$ are nonnegative rates.
For instance, every 
nucleotide $A$ mutates to $T$ at rate $v_W$, to $C$ at rate $v_S$, and
to $G$ at rate $w_S$.
The indices $W$ and $S$ refer to the classification of nucleotides
according to the strength of their link in double stranded DNA, the
link between $C$
and $G$ being strong ($S$) and the link between $A$ and $T$ being weak ($W$).

One recovers Tamura's matrix when $w_S\,v_W=v_S\,w_W$. On the other
hand, assuming that (b) holds, condition (a) corresponds to the
additional requirements that the two substitution rates from a purine
to $C$ coincide, and that the two substitution rates from a purine to
$T$ coincide.

As before, one can complete the matrix by diagonal elements, which represent fictitious rates of mutation from 
a nucleotide $x$ to $x$, and which leave the whole
process unchanged.
The full matrix is
$$
\left(
\begin{array}{cccc}
w_W & v_W & v_S & w_S
\\
v_{W} & w_W & w_S & v_S
\\
v_{W} & w_W & w_S & v_S
\\
w_W & v_{W} & v_S & w_S
\end{array}
\right).
$$ 
\subsection{Frequencies}
From now on, we assume that properties (a) and (b) hold, that the only nonzero YpR rates are 
$r_{S}$ and $r_{W}$ defined above, and we compute the frequencies of the nucleotides and of the YpR
dinucleotides.
\begin{definition}\label{d.sigma}
Introduce the parameters
$$
\ss_{S}:=v_{S}+w_{S},\quad \ss_{W}:=v_{W}+w_{W},\quad
\sv:=v_{S}+v_{W},\quad\sw:=w_{S}+w_{W},
$$
and
$$ \ss:=\ss_{S}+\ss_{W}=\sv+\sw.
$$
\end{definition}
Using the notations in section~\ref{s.tgc}, one gets
$$
u_{S}+u_{W}=\sv-\sw,\quad
v=2\sv,\quad
w=2\sw,\quad
s_{R}=s_{Y}=\ss.
$$
Using theorem~\ref{t.slin}, one gets
$$
M\times\left(\begin{array}{c}F(CG)\\
F(CA)\\
F(TA)\end{array}\right)=\frac12\,\left(\begin{array}{c}v_{S}\,\ss_{S}\\
v_{S}\,\ss_{W}+v_{W}\,\ss_{S}
\\
v_{W}\,\ss_{W}\end{array}\right),
$$
where the $3\times 3$ matrix $M$ is
$$
M:=\left(\begin{array}{ccc}
\ss\,(\ss+u_{S})+(\ss+v_{S})\,r_{W}
&
\ss\,u_{S}
&
-v_{S}\,r_{S}
\\
\ss\,u_{W}-(2v_S+\sw)\,r_{W}
&
\ss\,(3\sv+\sw)
&
\ss\,u_{S}-(2v_{W}+\sw)\,r_{S}
\\
-v_{W}\,r_{W}
&
\ss\,u_{W}
&
\ss\,(\ss+u_{W})+(\ss+v_{W})\,r_{S}
\end{array}
\right).$$
From there, tedious computations lead to the following formulas.
\begin{theorem}\label{t.wms}
In the symmetric R/Y + YpR model described by the parameters $v_S$,
$w_S$, $v_W$, $w_W$, $r_S$ and $r_W$, the frequencies of the YpR
dinucleotides at equilibrium
are
$$
F(xy)=D(xy)/(4D),
$$
where $D$ and $D(xy)$ are polynomial functions of the parameters
$(v_S,w_S,v_W,w_W,r_S,r_W)$,
homogeneous
of
degree $3$, and
defined as follows.
First,
$$
D:=D_{0}+r_{S}\,D_{S}+r_{W}\,D_{W}+r_{S}r_{W}\,D_{SW},
$$
with
$$
D_{0}:=\ss^{2}\,(\ss+2\sv),\quad
D_{SW}:=2\,(\ss+\sv),
$$
and
$$
D_{S}:=\ss\,(\ss+2\sv)+\ss\,w_{S}+\sv\,\ss_{S},\qquad
D_{W}:=\ss\,(\ss+2\sv)+\ss\,w_{W}+\sv\,\ss_{W}.
$$
As regards the dinucleotide $CG$,
$$
D(CG):=D_{0}(CG)+r_{S}\,D_{S}(CG),
$$
with
$$
D_{0}(CG):=(\ss+2\sv)\,\ss_{S}^{2},
\qquad
D_{S}(CG):=\ss\,w_{S}+\ss_{S}\,(\sv+2v_{S}).
$$
As regards the dinucleotide $TA$,
$$
D(TA):=D_{0}(TA)+r_{W}\,D_{W}(TA),
$$
with
$$
D_{0}(TA):=(\ss+2\sv)\,\ss_{W}^{2},
\qquad
D_{W}(TA):=\ss\,w_{W}+\ss_{W}\,(\sv+2v_{W}).
$$
As regards the dinucleotide $CA$ (and the dinucleotide $TG$),
$$
D(CA):=D_{0}(CA)+r_S\,D_S(CA)+r_{W}\,D_{W}(CA)+r_Sr_W\,D_{SW}(CA),
$$
with
$$
D_{0}(CA):=(\ss+2\sv)\,\ss_S\,\ss_{W},
\quad
D_{SW}(CA):=\ss+\sv,
$$
and
$$
D_{S}(CA):=\ss_S\,(\ss+2v_{W}),
\quad
D_{W}(CA):=\ss_W\,(\ss+2v_{S}).
$$
Finally,
$$
\ss\,F(C)=\ss_{S}/2-r_{W}\,F(CG)+r_{S}\,F(TA).
$$
\end{theorem}
One can check that, when $r_{S}=r_{W}=0$, $F(x)=F_{0}(x)$ and
$F(xy)=F_{0}(xy)$, where, for instance, 
$$
F_{0}(CG)=F_{0}(C)\,F_{0}(G),\quad F_{0}(C)=F_{0}(G)=\frac{\ss_{S}}{2\,\ss}.
$$
\subsection{$\cpgoe$ and $\tpaoe$ ratios}
We come back to the original model of YpR substitutions, where CpG is
the only active YpR dinucleotide, 
hence 
$$
r_{S}=0.
$$
In this case,
$$
F(CG)=\frac{D_{0}(CG)}{4\,(D_{0}+r_{W}\,D_{W})},
$$
and $F(C)$ is given by the relation
$$
\ss\,F(C)=\ss_{S}/2-r_{W}\,F(CG).
$$
This is enough to get some information about the ratio of the observed
and expected frequencies of CpG.
\begin{definition}
\label{d.oe}
Introduce
$$
\cpgoe:=\frac{F(CG)}{F(C)\,F(G)}.
$$
\end{definition}
\begin{prop}[$\cpgoe$ in a simple case]\label{p.cpgoe}
In the setting of theorem~\ref{t.wms}, assume furthermore that
$r_S=0$.
Then $\cpgoe\le1$ for every $r_W$, $\cpgoe$ is a non 
increasing function of $r_W$,
$\cpgoe\to1$ when $r_W\to0$, and $\cpgoe\to0$ when $r_W\to\infty$. 
Furthermore, when $r_W\to0$,
$$
\cpgoe=1-r_{W}\,K_{CG}+o(r_{W}),
$$
where $K_{CG}$ is positive and defined as
$$
K_{CG}:=\frac{(\ss+3\sv)\,\ss_{W}+\ss\,w_{W}}{\ss^{2}(\ss+2\sv)}.
$$
\end{prop}
Another quantity of interest is the ratio of the observed
and expected frequencies of TpA.
\begin{definition}
\label{d.oe2}
Introduce
$$
\tpaoe:=\frac{F(TA)}{F(T)\,F(A)}.
$$
\end{definition}
In the setting of theorem~\ref{t.wms}, and even if one assumes furthermore that
$r_S=0$, the situation is less clear than for $\cpgoe$.
For instance, one can show that, when $r_W\to0$,
$$
\tpaoe=1-r_{W}\,K_{TA}+o(r_{W}),
$$
where $K_{TA}$ is defined as
$$K_{TA}:=\frac{L_{TA}}{\ss^{2}\,\sw^2\,(\ss+2\sv)},
$$ 
and
\begin{eqnarray*}
L_{TA} & := &
\ss_W\ss_S^2(\ss+2\sv)+\ss_W^2(\ss(\ss+2\sv)+\ss
w_W+v\ss_W)
\\
& & \hspace{20ex}-\ss^2(\ss w_W+v\ss_W)-2\ss^2\ss_W v_W.
\end{eqnarray*}
The sign of $L_{TA}$ is difficult to decipher, in fact assume that
there exists $c$ such that
$$
w_S=c\,v_S,\quad w_W=c\,v_W.
$$
Then the expression of $L_{TA}$ reduces to 
$$
L_{TA}=\ss\,\ss_W^2\,((3+c)\,v_W-(1+c)\,v_S).
$$
This shows the following result.
\begin{prop}[Values of $\tpaoe$]\label{p.tpa}
In the general setting of proposition~\ref{p.cpgoe}, both cases
$\tpaoe\le1$ and
$\tpaoe\ge1$ are possible. However, when $w_S=w_W$ and $v_S=v_W$, 
$$\tpaoe\le1.$$
\end{prop}

\section{The simplest model}\label{s.soa}

In this section, we provide the values of the $16$ frequencies of
dinucleotides at equilibrium. To avoid awkward formulas, we consider
the simple non trivial case.
\begin{definition}\label{d.mss}
The simplest R/Y + YpR model is such that,
 for every nucleotide $x$, 
$$
v_{x}=w_{x}=1,
$$ 
and such that all YpR substitutions but those starting from CpG are excluded, i.e.
$$
r^C_{G}=r^{A}_{T}=r^G_{C}=r^T_{A}=r^{T}_{G}=r^{A}_{C}=0,
$$
and such that the rates of CpG to CpA and CpG to TpG are equal, that is,
$$r^{C}_{A} = r^{G}_{T}=:\varrho.$$ 
\end{definition}
In this section, we consider the simplest model.

\subsection{On symbolic resolutions}
Theoretically, one has to solve an
appropriate linear system related to the dynamics on the discrete
circle with $N+2=4$ vertices, that is, of size $4^{N+2}=256$.  The
translation invariance yields a system of size $m(N+2)=70$, see
section~\ref{ss.pf}.
Using the invariance with respect to both
translations of the discrete circle and nucleotide complementarity,
the size of the linear system to be solved is further 
reduced to $42$. (This is because $14$ of the $70$ classes that the
invariance by translations induces, are invariant by the 
nucleotide complementarity
as well, the $56$ other classes being grouped into pairs. 
We omit the details of this enumeration.)

We solved this $42\times42$ 
system symbolically, using Maple\texttrademark.  We
computed the full invariant distribution of $(X^{\II}(s))_{s \geq 0}$
with $\II:=\{1,2,3,4\}$ but we only give the equilibrium
frequencies of dinucleotides because these are the quantities of
greatest interest. Checking the results by human computations seemed
prohibitively time-consuming and tedious but, to confirm the
validity of the formulas, we performed some tests. In particular, we
compared the formulas with the following.
\begin{itemize}
\item The exact formulas obtained by human computations for YpR dinucleotides.
\item The results obtained by numerically solving the system with
  MATLAB\textregistered\ for various settings of the parameters.
\item The results of extensive Monte-Carlo simulations, usually with
  $10^{8}$ runs for each setting of the parameters. 
\end{itemize}
All these tests confirmed the values given below.

\subsection{Frequencies}

We first recall the values of the frequencies of nucleotides, deduced
from previous sections.
\begin{definition}
\label{d.kk}
Introduce
$$
K_0(xy)=4\,U(xy)+2\,R(x)+Y(x)+R(y)+2\,Y(y).
$$
where
$$
U=\un_{TG}+\un_{CA}-2\,\un_{CG},
\quad
R=\un_A-\un_G,
\quad
Y=\un_T-\un_C.
$$
Introduce
$$
a(\varrho):=\frac3{96+19\varrho},
\quad
b(\varrho):=\frac4{32+10\varrho}.
$$
\end{definition}
\begin{definition}
For every polynucleotide $x_1\cdots x_k$, define a function
$K(x_1\cdots x_k)$ by the relation
$$
F(x_1\cdots x_k)=:\frac1{4^k}\,
\left(1+\varrho\,\frac{K(x_1\cdots x_k)}{32+10\varrho}\right).
$$
\end{definition}
\begin{theorem}\label{tk0}
{\bf (1)}
For every nucleotide $x$, $K(x)$ does not depend on $\varrho$, and
$$
K(x)=2\,(R(x)+Y(x)).
$$
Hence $K(A)=K(T)=2$ and $K(C)=K(G)=-2$.
\\
{\bf (2)}
For every YpR dinucleotide $xy$, $K(xy)$ does not depend on $\varrho$,
and $K(xy)=K_0(xy)$. Hence
$$
K(CG)=-10,\quad K(CA)=K(TG)=4,\quad K(TA)=2.
$$
{\bf (3)}
For every dinucleotide, $K(xy)\to K_0(xy)$ when
$\varrho\to0$. 
\end{theorem}
Part (3) reads as
$$
K_0(GG)=K_0(CC)=-3,
\quad
K_0(TT)=K_0(AA)=3,
$$
$$
K_0(AG)=K_0(CT)=1,
\quad
K_0(TC)=K_0(GA)=-1,
$$
and
$$
K_0(AC)=K_0(GT)=0,\quad
K_0(AT)=4,\quad
K_0(GC)=-4.
$$
Here is a consequence of theorem~\ref{tk0}.
\begin{prop}
The nucleotides $C$ and $G$ are always less frequent than $A$ and
$T$. More precisely, for every positive $\varrho$,
$$
20 \%\le P(C)=P(G)<25\%<P(A)=P(T)\le30\%.
$$ 
Furthermore, the dinucleotides
$CG$ and $TA$ are repulsive and the dinucleotides $CA$ and $TG$ are
attractive, in the sense that
$$
F(CG)\le F(C)\,F(G),\quad F(TA)\le F(T)\,F(A),
$$
and
$$
F(CA)\ge F(C)\,F(A),\quad F(TG)\ge F(T)\,F(G).
$$
\end{prop}
\begin{definition}\label{d.tk2}
For every dinucleotide $xy$,  introduce
$K_1(xy)$ as
$$
K(xy)=:K_0(xy)+\varrho\,K_1(xy).
$$
\end{definition}
\begin{theorem}
\label{tk1}
For every dinucleotide $xy$,
$K_1(xy)$ assumes one of the five values $0$, $\pm a(\varrho)$, and
$\pm b(\varrho)$.
More precisely,
$$
K_1(xy)=a(\varrho)\,(R(x)\,R(y)+Y(x)\,Y(y))+b(\varrho)\,R(x)\,Y(y).
$$
\end{theorem}
As regards, for instance, the dinucleotides $Ax$, this means that
$$
K_1(AA)=a(\varrho),
\
K_1(AC)=-b(\varrho),
\
K_1(AG)=-a(\varrho),
\
K_1(AT)=b(\varrho).
$$
Going back to frequencies, this reads as
\begin{eqnarray*}
F(AA) & = & \frac1{16}\,\left(1+\frac{\varrho}{32+10\varrho}\,\left(3+
\frac{3\varrho}{96+19\varrho}\right)\right),
\\
F(AC) & = & \frac1{16}\,\left(1+\frac{\varrho}{32+10\varrho}\,\left(0-\frac{4\varrho}{32+10\varrho}\right)\right),
\\
F(AG) & = & \frac1{16}\,\left(1+\frac{\varrho}{32+10\varrho}\,\left(1-\frac{3\varrho}{96+19\varrho}\right)\right),
\\
F(AT) & = &
\frac1{16}\,\left(1+\frac{\varrho}{32+10\varrho}\,\left(4+\frac{4\varrho}{32+10\varrho}\right)\right).
\end{eqnarray*}
Similar formulas are available for the $12$ other dinucleotides.

\subsection{Remarks}

One sees that, when $\varrho\to\infty$, $\varrho\,K_1(xy)$ converges
to a nondegenerate limit.
\begin{corollary}\label{c.toi}
When $\varrho\to\infty$, 
$K(xy)$ converges to
$K_{\infty}(xy)$, with
$$
K_{\infty}(xy):=K_0(xy)+\frac3{19}\,\big(R(x)\,R(y)+Y(x)\,Y(y)\big)
+\frac25\,R(x)\,Y(y).
$$
\end{corollary}
This means that
$$
F(xy)\to 
F_{\infty}(xy):=\frac1{16}\,\left(1+\frac{K_{\infty}(xy)}{32}\right).
$$
For instance, 
$$F(C)=F(G)\to\frac15,
\quad
F(A)=F(T)\to\frac3{10},
$$ 
and
$$
F(CG)\to0,\qquad
F(CA)=F(TG)\to\frac7{80},\qquad
F(TA)\to\frac3{40}.
$$
When $\varrho=0$, one recovers the i.i.d.\ evolution, hence
$F(x_1\cdots x_k)=1/4^{k}$
for every polynucleotide $x_1\cdots x_k$.

The model makes sense for every $\varrho\ge-1$. For instance,
$\varrho=-1$ forbids the transitions from $CG$ to $TG$ or $CA$, but
allows the transitions from $CG$ to $AG$, $GG$, $CC$ and
$CT$. Although the model does not make sense when $\varrho<-1$, the
frequencies that we are able to compute are all positive, even
formally for
some values of $\varrho<-1$. It may happen that, as soon as
$\varrho<-1$, the expressions of the frequencies of some
polynucleotides are indeed negative.
In any case, when $\varrho=-1$,
$$
F(C)=F(G)=\frac5{22},\quad
F(A)=F(T)=\frac6{22}.
$$

\section{Continuous dynamics}
\label{s.cd}
In this section, we show that one can study the evolution of the
YpR frequencies, using essentially the techniques used to describe the
 statics. 
Consider for instance the simplest model, see definition~\ref{d.mss}.
For every YpR dinucleotide $xy$, the frequency $F(xy)(s)$ of $xy$ at
 time $s$ satisfies
\begin{eqnarray*}
\frac\dd{\dd
  s}\,F(xy)(s) & = & F(x)(s)+F(y)(s)-8\,F(xy)(s)+
\\
& &
  {}+\varrho\,F(CG)(s)\,(\un_{TG}+\un_{CA})(xy)
-2\varrho\,F(CG)(s)\,\un_{CG}(xy).   
\end{eqnarray*}
Likewise, for every nucleotide $x$,
\begin{equation}
\label{e3}
\frac\dd{\dd s}\,F(x)(s)=-4\,F(x)(s)+1+\eps(x)\,\varrho\,F(CG)(s),
\end{equation}
where $\eps(C)=\eps(G)=+1$ and $\eps(A)=\eps(T)=-1$. 
Hence, $F(CG)(s)$ satisfies the second order evolution equation
$$
\frac{\dd^{2}}{\dd s^{2}}\,F(CG)(s)+2\,(\varrho+6)\,\frac\dd{\dd s}\,F(CG)(s)
+2\,(16+5\varrho)\,F(CG)(s)=2.
$$
One can check that the two roots of the characteristic polynomial of this linear differential equation
have negative real parts. Hence
$F(CG)(s)$ indeed converges when $s\to\infty$
to the fixed point $1/(16+5\varrho)$ of the equation, that is, to the
stationary value $F(CG)$.
Plugging this into (\ref{e3}) for every value of $x$ yields the convergence of 
$4\,F(x)(s)$, when $s\to\infty$, to the value
$$
1+\eps(x)\,\varrho\,F(CG)=4\,F(x).
$$

\section*{Part C Simulation}
\addcontentsline{toc}{section}{Part C Simulation}
\section{Coupling from the past}
\label{s:CFTP}

A consequence of the construction of the previous sections is that we can simulate the restriction of the dynamics on $\Z$ to any finite interval of sites of length $n$, 
without truncation errors due to neglecting the influence of remote sites.
 One adds a site to the
left and a site to the right, one performs simulations for the
system on these $n+2$ sites, and the projection on the $n$ original
sites yields the desired simulation.
 
In this section, we show how the coupling from the past (CFTP)
methodology of Propp and Wilson~\cite{ProWil} can be applied in our
context. Our motivation is two-fold. First, estimates about the
coupling times automatically yield estimates on the speed of
convergence of the dynamics to the stationary distribution. In our
context, this applies to the speed of convergence of
$\prt^{a,b}(X_{s})$ to $\prt^{a,b}(\mu)$.  Second, the CFTP technique
allows to sample exactly from $\prt^{a,b}(\mu)$. Despite the results
of the previous sections, which show that the obtention of exact
expressions of $\prt^{a,b}(\mu)$ amounts to the inversion of a linear
system, this task becomes computationally infeasible as soon as the
number $b-a+1$ of sites is large, say greater than $6$. Hence,
Monte-Carlo simulations are still useful, if only to confirm the
results obtained by inverting the linear system!

In the whole section, we fix $\II=\ab$.

\subsection{Coupling events}\label{ss.ce}

We first define the notions of coupling events and locked sites.
\begin{definition}[Coupling events]\label{d.ce} 
We say that a coupling event occurs at site $i$ and times $(s_{1}, s_{2},
s_{3})$ if the following assertions hold.
\begin{itemize}
\item $s_{1} > s_{2}$ and $s_{3} > s_{2} $,
\item $-s_{1} $ belongs to $\UU^{A}_{\lef(i)}  \cup \UU^{G}_{\lef(i)} \cup \VV^{A}_{\lef(i)}  \cup
\VV^{G}_{\lef(i)}$,
\item $-s_{3} $ belongs to $  \UU^{C}_{\rig(i)}  \cup \UU^{T}_{\rig(i)}  \cup  \VV^{C}_{\rig(i)}  \cup
\VV^{T}_{\rig(i)} $,
\item  $-s_{2} $ belongs to $\bigcup_{z \in \A}  \UU^{z}_{i}$,
\item $(-s_{1},-s_{2}) \cap    \bigcup_{z \in \A }  \UU^{z}_{\lef(i)}\cup\VV^{z}_{\lef(i)} $ is empty,
\item $(-s_{3},-s_{2}) \cap    \bigcup_{z \in \A }  \UU^{z}_{\rig(i)}\cup\VV^{z}_{\rig(i)} $ is empty.
\end{itemize}
\end{definition}
\begin{definition}[Locked sites]\label{d.ls}
We say that the site $i$ is locked at times $(u,v)$ if, for every times $\sigma$ and
$s$ such that $\sigma \le -u$ and $s \ge  -v$, 
the set
$$
\Phi_{\sigma} (\A^{\II},\xi,i,s)
$$
contains but one single element. In other words, $\Phi_{\sigma}(x,\xi,i,s)=\Phi_{\sigma}(x',\xi,i,s)$ 
for every initial configurations $x$ and $x'$ in $\A^{\II}$.
\end{definition}
Recall that $\Phi_\sigma$ is introduced in section~\ref{s.tpoaz}.
Definitions~\ref{d.ce} and \ref{d.ls} both involve $\II$, 
through the definitions of $\lef(\cdot)$ and 
$\rig(\cdot)$. 
Our next proposition shows that the notions of coupling event and locked
site are related. The proofs of the results in section~\ref{s:CFTP}
are in section~\ref{ss.12p}.
\begin{prop}\label{p:coa}
If a coupling event occurs at site $i$ and times  $(s_{1}, s_{2}, s_{3})$, then
$i$ is locked at times  $(s_{4}, s_{2})$, with $s_{4}:=\max(s_{1},s_{3})$.   
\end{prop}
In words, for every initial condition imposed before time $-\max(s_{1},s_{3})$, the $i$th 
coordinate is the same after
time $-s_{2}$. A consequence is that all the trajectories of the process have coalesced as far as site $i$ is
concerned.

\subsection{Locking times}
\label{ss.lt}
We wish to estimate, or at least to control, the time that is needed to lock a given collection of sites.
We start with one site.
\begin{definition}[Locking times]\label{d.lt}
Let $T_{i}$  be defined  as the supremum of the times $s_{4}$ such that there exists $s_{1}$, $s_{2}$ 
and $s_{3}$, with the following properties:
$s_{2}\ge0$, $s_{4}=\max(s_{1},s_{3})$, and a
coupling event occurs at site $i$ and at times $(s_{1}, s_{2}, s_{3})$.
\end{definition}
We now define additional numerical parameters.
\begin{definition}[Combined rates]\label{d.cr}
Introduce
$$
\vara_{R}:=\sum_{z = A,G} \min(c_{z},v_{z}),\quad
\vara_{Y}:=\sum_{z = T,C} \min(c_{z},v_{z}),\quad
\vara:=\vara_{R}+\vara_{Y},
$$
and
$$
\varb_{R}:=\sum_{z = A,G} (v_{z}-c_{z})^{+},\quad
\varb_{Y}:=\sum_{z = T,C} (v_{z}-c_{z})^{+},\quad
\varb:=\varb_{R}+\varb_{Y}.
$$
\end{definition}
Note that $\vara_{R}+ \varb_{R}=v_{A}+v_{G}$ and $\vara_{Y}+\varb_{Y}=v_{C}+v_{T}$.
\begin{definition}\label{d.bata}
Let 
$$
\bata:=t_{Y}\,t_{R}=  
\frac{(v_{A}+v_{G}) \times (v_{C}+v_{T}) }{ (v_{A}+v_{G}+v_{C}+v_{T} )^{2}}.
$$ 
\end{definition}
These parameters allow us to control the distribution of the locking times, as follows.
\begin{prop}\label{p:sadom}
Each locking time $T_i$ is stochastically dominated by 
the random variable 
$$
H_{1}+\cdots+H_{Z},
$$
where $(H_{k})_{k\ge1}$ is a sequence of i.i.d.\ Gamma $(3,\vara)$ random variables, 
and $Z\ge1$ is a geometric random variable with parameter $\bata$, independent from $(H_{k})_{k\ge1}$.
In particular, $T_{i}$ is almost surely finite and $\E(T_{i})\le3/(\vara\bata)$.
\end{prop}

We now proceed to bound the tail of $T_{i}$. For each $k\ge1$, the distribution of 
$H_{1}+\ldots+H_{k}$ is Gamma $(3k,\vara)$.

\begin{definition}\label{d.nalpha}
Introduce $n_{\bata}=-4\log(1-\bata)/\bata$, hence $n_{\bata}$ is finite and $n_{\bata}\ge4$.
\end{definition}
\begin{lemma}\label{l:queue}
For every nonnegative integer $N$,
$$
\P (H_{1}+\cdots+H_{Z} \ge Nn_{\bata}/\vara) \le  2\,(1-\bata)^{N}.
$$
\end{lemma}

Our interest lies in the coupling time of whole intervals,
defined below.
\begin{definition}[Locking times of intervals]\label{d.cti}
The locking time $T_{a,b}$ of the sites in the interior of
the interval $\II=\{a,\ldots,b\}$ is
$$
T_{a,b} = \max\{T_{i}\,;\,a+1\le i\le b-1\}.
$$
\end{definition}
\begin{definition}\label{d.kab}
Introduce the integer 
$$
k_{a,b}=\pentsup{\frac{b-a-1}{3}}.
$$
\end{definition}
\begin{prop}\label{p.lti}
For every integer $N$,
$$\Q ( T_{a,b}   \ge  N n_{\bata} /\vara) \le   3   \left[ 1    -       \left( 1    - 
2\,(1-\bata)^{N}  \right)^{k_{a,b}} \right].
$$
\end{prop}
\begin{remark}
The bound above has the nice property that it does not involve the rates $r^y_{z}$
of the YpR mutations except, through the $c_{z}$, when they are negative.  
\end{remark}

More readable forms of proposition~\ref{p.lti} might be proposition~\ref{p.rea}
and corollary~\ref{c.rea} below.
\begin{definition}\label{d.tab}
Let $$t_{a,b} =4\log(6\,k_{a,b}).$$
\end{definition}
\begin{prop}[Control upon the locking times]\label{p.rea}
For every $t$, 
$$\Q((\bata\vara)\, T_{a,b}   \ge  t_{a,b}-\log(1-\bata)+t ) \le \exp(-t/4).$$
\end{prop}
\begin{definition}\label{d.lti}
Let $T_{(n)}$ denote the locking time of $n$ consecutive sites.
\end{definition}
For instance $T_{(n)}$ is distributed 
as $T_{0,n+1}$. Note that $3k_{0,n+1}\le n+2$ and that $\bata\le\frac14$, hence $4\log2-\log(1-\alpha)\le6\log2$.
This yields the following corollary.
\begin{corollary}\label{c.rea}
For every $t$, 
$$\Q((\bata\vara)\, T_{(n)}   \ge  \log(n+2)+6\log(2)+t ) \le \exp(-t/4).$$
\end{corollary}

\subsection{Consequences}

The now traditional Propp-Wilson method induces that, whatever the initial condition $\bxx$ in $\A^{\II}$ with
$\II=\{a,  \ldots, b \}$,
for every $J\subset\II$, 
$$\left[\varphi^{a,b}_{-T } (\bxx,\xi, j,0) \right]_{j \in J},\quad
\mbox{where}\ T=\max_{j \in J} T_{j},
$$ 
is distributed according to the
projection of the 
stationary distribution $\mu_{\II}$ on $\A^{J}$. In particular,
the distribution of 
$$
\prt^{a+1,b-1} \left(\varphi^{a,b}_{-T_{a,b}} (\bxx,\xi)(0)\right)
$$ 
is 
$  \prt^{a+1,b-1} (\mu_{a,b})$, that is,  by proposition~\ref{p.r}, $\prt^{a+1,b-1} (\mu)$.
We state this as a proposition. 
\begin{prop} For every $a \le b$,
the distribution of 
$
\prt^{a,b} \left(\varphi^{a-1,b+1}_{-T_{a-1,b+1}} (\bxx,\xi)(0)\right)
$ is 
$\prt^{a,b} (\mu)$.
\end{prop}
Hence proposition~\ref{p.rea} yields the result below.
\begin{prop}
For every positive $t$, the distance in total variation between the distribution
$\prt^{a,b}(X^{\bxx}_{s})$ at time
$$s=(t_{a-1,b+1}-\log(1-\bata)+t)/(\bata\vara),
$$
and the limiting distribution $\prt^{a,b}(\mu)$,  is at most $\exp(-t/4)$.

\end{prop}

\subsection{Proofs}
\label{ss.12p}
\begin{proof}[\ of proposition~\ref{p:coa}]
We shall in fact prove the following assertion: assume that site $i$ is locked at times $(s_{1}, s_{2},
s_{3})$, then, for every $\sigma \le -s_{4}$ and $s \ge  -s_{2}$, the sets
$$
\varrho(\Phi_{\sigma} (\A^{\II},\xi,\lef(i),s)), \quad \Phi_{\sigma} (\A^{\II},\xi,i,s), \quad
\eta(\Phi_{\sigma}
(\A^{\II},\xi,\rig(i),s))$$
are all singletons.

We first check the claim when $s=-s_{2}$. By the definition of $s_{1}$, at time $-s_{1}$, either a move of 
type $\typeu$ occurs at
 site $\lef(i)$, yielding an $A$ or a $G$ unconditionally, or a move of type $\typev$ occurs, namely a transversion to 
a purine, yielding a purine if site $i$ was not already occupied by a purine. As a consequence,  the set
$$
\varrho \left( \Phi_{\sigma} (\A^{\II},\cdot,\lef(i),-s_{1}) \right)
$$ 
is a singleton. Once again by the
definitions, we ruled out the possibility that any move of type $\typeu$ or $\typev$ occurred at site 
$\lef(i)$
between the times $-s_{1}$ and $-s_{2}$. Furthermore, moves of type
$\typew$, $\typer$ and $\typeq$, when applied to a purine, can only
yield a (possibly different) purine. This implies that the set
$$
\varrho(\Phi_{\sigma} (\A^{\II},\xi,\lef(i),-s_{2}))
$$
is a singleton as well. 
The same argument applies symmetrically to $\eta(\Phi_{\sigma} (\A^{\II},\xi,\rig(i),-s_{2}))$.
 
As regards $\Phi_{\sigma} (\A^{\II},\xi,i,-s_{2})$, this is a singleton
since a move of type $\typeu$ occurs at site $i$ at time $-s_{2}$. Lemma~\ref{l.d1}
above shows that, for every $x$ in $\A^{\II}$, 
the values of 
$$
\varrho(\Phi_{\sigma} (\bxx,\xi,\lef(i),s)),
\quad
\Phi_{\sigma} (\bxx,\xi,i,s),
\quad
\eta(\Phi_{\sigma} (\bxx,\xi,\rig(i),s)),
$$ 
for every $s \ge 
-s_{2}$, are completely determined by $\xi$ and by the values of
$$
\varrho(\Phi_{\sigma} (\bxx,\xi,\lef(i),-s_{2})),
\quad
\Phi_{\sigma} (\bxx,\xi,i,-s_{2}),
\quad
\eta(\Phi_{\sigma} (\bxx,\xi,\rig(i),-s_{2})). 
$$ 
Since these values are the same for every $\bxx$  in $\A^{\II}$, so is the case for
 $\varrho(\Phi_{\sigma} (\cdot,\xi,\lef(i),s))$, $\Phi_{\sigma} (\cdot,\xi,i,s)$, and 
$\eta(\Phi_{\sigma}
(\cdot,\xi,\rig(i),s))$, for every $s \ge  -s_{2}$.
This concludes the proof.
\end{proof}
\begin{proof}[\ of proposition~\ref{p:sadom}]
Recall the convention that $\sup \emptyset = - \infty$.
Define $M_{0}:=0$, and, inductively for $k \ge 1$,
\begin{itemize}
\item $-L_{k} := \sup \ (-\infty,-M_{k-1}) \cap \left(  \bigcup_{z \in \A }  \UU^{z}_{i}     \right)$,
\item $-U_{k} := \sup \ (-\infty, -L_{k}) \cap    \left( \bigcup_{z \in \A}  \UU^{z}_{\lef(i)}\cup\VV^{z}_{\lef(i)}
 \right)$, 
\item $-V_{k} := \sup \ (-\infty, -L_{k}) \cap    \left( \bigcup_{z \in \A}  \UU^{z}_{\rig(i)}\cup\VV^{z}_{\rig(i)}
 \right)$, 
\item $-M_{k} := -\max (U_{k},V_{k})$.
\end{itemize}
Define $K$ as the smallest integer $k \ge 1$ such that
$$-U_{k}\ \mbox{belongs to}\  \UU^{A}_{\lef(i)}  \cup    \UU^{G}_{\lef(i)}     \cup  \VV^{A}_{\lef(i)}  \cup   
\VV^{G}_{\lef(i)},
$$     
and 
$$
-V_{k} \ \mbox{belongs to}\  \UU^{C}_{\rig(i)}  \cup  \UU^{T}_{\rig(i)}  \cup   \VV^{C}_{\rig(i)}  \cup 
\VV^{T}_{\rig(i)}.
$$ 
Then, provided that $K$ is finite,  a coupling event occurs at site $i$ and times $(U_{K},L_{K},V_{K})$, hence
 $T_{i} \le M_{K}$ as soon as $K$ is finite. 
Furthermore, standard properties of Poisson processes and the independence of the Poisson processes that are
associated to different sites show that the sequence
$$\big(L_{k}-M_{k-1},  
U_{k} - L_{k} , V_{k} - L_{k}\big)_{k\ge 1}
$$ 
is i.i.d., and that,
for every given $k\ge1$,  $L_{k}-M_{k-1}$,  
$U_{k} - L_{k}$ and $V_{k} - L_{k}$ are mutually independent and
exponentially distributed with parameters $\vara$,
 $\vara+\varb$, and $\vara+\varb$ respectively. 
Finally, $K\ge1$
is independent from  $(L_{k}-M_{k-1}, U_{k} - L_{k} , V_{k} -
L_{k})_{k\ge 1}$, and geometrically distributed with 
 parameter $\bata$ and expectation $1/\bata$.
Writing 
$$
M_{k}-M_{k-1} = L_{k} - M_{k-1}    + \max(U_{k} - L_{k},
V_{k} - L_{k}),
$$ 
and recalling that  $T_{i} \le M_{K}$, one gets
$$
\E(T_i)\le\frac1{\bata}\,\left(\frac1{\vara}+\frac3{2(\vara+\varb)}\right).
$$
Simpler upper bounds obtain as follows.
Since $  \vara+\varb   \ge \vara$, the distribution of the random
variable
$\max( U_{k} - L_{k}, V_{k} -
L_{k})$ is (crudely) dominated by the distribution of the sum of two independent $\vara$ exponential
random
variables, hence the distribution of $T_i$ is dominated by the distribution of the sum of three independent $\vara$ exponential
random
variables.
\end{proof}
One sees that
$$
\E(T_i)\le\frac{5}{2\bata\vara}.
$$
\begin{proof}[\ of lemma~\ref{l:queue}]
By the homogeneity of the Gamma distributions, we assume that
$\vara=1$.
The nonnegativity of the random variables $(H_{k})_{k}$ implies that,
for every integer $k\ge0$ and every real number $t\ge0$, 
$$
\P(H_{1}+\cdots+H_{Z}\ge t)\le\P(Z\ge k+1)+\P(H_{1}+\cdots+H_{k}\ge t).
$$
The first term on the right hand side is $(1-\bata)^{k}$.
By Cram\'er's bound and the value of the Laplace transform of the standard
exponential distribution, evaluated at $0\le
u<1$,
 the
second term is at most 
$$
(1-u)^{-3k}\,\exp(-ut).
$$
Assume that $t=n_{\bata}N$ for an integer $N$, and choose $u=\bata$ and $k=N$.
Then the proof is complete, since for these values,
$$
(1-\bata)^{k}=(1-u)^{-3k}\,\exp(-ut)=(1-\bata)^{N}.$$
\end{proof}
\begin{proof}[\ of proposition~\ref{p.lti}]
Write $\II = I_{0} \cup I_{1} \cup I_{2}$,
where $I_{j}$ collects the sites in $\II$ that are equal to $j$ modulo 3.
For $j=0$, $1$ and $2$, let $T_{(j)}=\max\{T_{i}\,;\,i\in I_{j}\}$.
 
By the independence properties of the collection $(T_{i})_{i}$, for each $j$, the random variables
$(T_{i})_{i\in I_{j}}$ are i.i.d. Furthermore, for each $j$,
$T_{(j)}$ involves at most $k_{a,b}$ sites in $\II$.
Hence, for every nonnegative $t$,
$$
\Q(T_{(j)}\ge t)\le1-(1-\Q(T_{i}\ge t))^{k_{a,b}}.
$$
Since $T_{a,b}$ is the maximum of 
the three random variables $T_{(j)}$, 
$$
\Q(T_{a,b}\ge t)\le \Q(T_{(0)}\ge t)+\Q(T_{(1)}\ge t)+\Q(T_{(2)}\ge t).
$$
One concludes, using the upper bound of $\Q(T_{i}\ge t)$ in lemma~\ref{l:queue} above.
\end{proof}

\section{Practical issues}\label{s:bepractical}

This section is devoted to some practical issues related to an
effective implementation of the CFTP method of simulation of R/Y+YpR
systems. We consider two slightly different versions of the method and
we give short presentations of both. The key point of each version is
that an efficient detection of the coalescence is at hand. For the
sake of readability, we do not provide detailed pseudo-codes but only
the basic schemes used to implement the methods.

\subsection{Additional notations}
\label{ss.adn}
Fix $\xi$ in $\Omega_{1}$ and $x$ in $\A^{\II }$ with $\II=\{a,\ldots,b\}$.
Let
$$\T' =  (-\infty,0)\cap\bigcup_{i\in\II} \bigcup_{z \in \A} \UU^{z}_{i}(\xi)\cup\VV^{z}_{i}(\xi),$$
and
$$\T'' =  (-\infty,0)\cap\bigcup_{i\in\II} \bigcup_{z \in \A} \WW^{z}_{i}(\xi)\cup\RR^{z}_{i}(\xi)\cup\QQ^{z}_{i}(\xi).$$
Let $t'_{-1}=0$ and  $(-t'_{0} , -t'_{1},\cdots)$ denote the ordered list of points in $\T'$, that is:
$$ 
\T'=\{-t'_0,-t'_1,-t'_2,\cdots\}\quad\mbox{with}\quad
0 > -t'_{0} > -t'_{1} > \cdots.
$$ For every $n \ge 0$, we describe the move that occurs at time
$t'_{n}$ through the site $c'_{n}$ where the move occurs and through
the description $M'_{n}=(z'_{n},f'_{n})$ of the move, as defined above
in sections~\ref{ss.cotp} and ~\ref{ss.prop}.

For every $n \ge 0,$ let $N_{n}$ denote the number of points of $\T''$
in the interval $(-t'_{n},-t'_{n-1})$.  When $N_{n}\ge1$, let
$$\T'' \cap (-t'_{n},-t'_{n-1}) = \left\{ -t''_{n,1} < \ldots <  -t''_{n, N_{n}} \right\}.
$$
Finally, we describe the move which occurs at time $t''_{n,k}$ through the site $c''_{n,k}$ where the move
occurs and through the description
$M''_{n,k}=(z''_{n,k},f''_{n,k})$ of the move, as above.

\subsection{First algorithm}
\label{ss.fa}
We describe a routine which yields a random element of $\A^{ \{a+1,\ldots,b-1\} }$ with distribution
$\prt^{a+1,b-1}(\mu)$.
\begin{tt}
\begin{flushleft}
*** Coalescence detection ***\\
Let $n:=0$; \\
Until for every $i$ in $\{a+1,\ldots,b-1\}$, $T_{i} \le t'_{n-1}$,
do: 
\\
\hspace{3ex}
$\{$ Generate and store $M'_{n}$ and
$c'_{n}$; Let    $n:=n+1$;     $\}$; \\
***  Sampling ***\\
 Let $x := (A,\ldots,A)$; \\
 For $k$  going backwards from $n-1$ to $0$ do:  \\
\hspace{3ex} $\{$ Let $x := \gamma(x, M'_{k}, c'_{k})$; \\
\hspace{3ex} Generate $N_{k}$; \\
 \hspace{3ex} For $m$ going from $1$ to $N_{k}$ do: 
\\
\hspace{6ex} $\{$ Generate
 $M''_{k,m}$ and $c''_{k,m}$ ; 
Let $x:=  \gamma(x, M''_{k,m}, c''_{k,m})$; $\}$    $\}$ 
\\
  
Return $\prt^{a+1,b-1}(x)$.
 
 \end{flushleft}
 \end{tt}

 The feasibility of the above routine relies on several facts. First,
  it is easy to generate realizations of the random variables $c'_{n}$
  and $M'_{n}$, whose distributions are explicitly known and standard.
  Second, one can check whether or not $T_{i} \le t'_{n}$, knowing
  only the sequence $c'_{k}$ and $M'_{k}$ for $k$ between $0$ and $n$,
  and this can be done in a step-by-step way, updating information
  about sites as time goes backwards and new moves are introduced.

 One advantage of this method is that one does not have to generate
 the random variables $M''_{k,m}, c''_{k,m}$ on the coalescence
 detection pass nor to store them, but just to compute their effect on
 $x$ in a step-by-step way during the sampling pass. This helps
 keeping memory storage requirements and execution time to a minimum.

\subsection{Second algorithm}
\label{ss.sa}
We need some additional definitions, because the detection of coalescence
devised in this second algorithm uses a slightly different technique
from the one in the first algorithm.
\begin{definition}[Coupling event on an interval]\label{d.ceoai}
Let $J \subset \II$. 
A coupling event of type $(s_{1},s_{2},s_{3},J)$
occurs at site $i$ when the following conditions hold:
\begin{itemize}
\item $s_{1} > s_{2}$ and $s_{3} > s_{2} $,
\item $-s_{1} \in  \UU^{A}_{\lef(i)}  \cup \UU^{G}_{\lef(i)} \cup \VV^{A}_{\lef(i)}  \cup
\VV^{G}_{\lef(i)}$ or $\lef(i) \in J$,
\item $-s_{3} \in  \UU^{C}_{\rig(i)}  \cup \UU^{T}_{\rig(i)}  \cup  \VV^{C}_{\rig(i)}  \cup
\VV^{T}_{\rig(i)} $  or $\rig(i) \in J$,
\item  $-s_{2}  \in \bigcup_{z \in \A}  \UU^{z}_{i}$,
\item $(-s_{1},-s_{2}) \cap    \bigcup_{z \in \A}  \UU^{z}_{\lef(i)}\cup\VV^{z}_{\lef(i)} $ is empty
or $\lef(i) \in J$,
\item $(-s_{3},-s_{2}) \cap    \bigcup_{z \in \A}  \UU^{z}_{\rig(i)}\cap\VV^{z}_{\rig(i)} $ is empty
or $\rig(i)\in J$.
\end{itemize}
\end{definition}
The key observation is the following: assume that a coupling event of type $(s_{1},s_{2},s_{3},J)$ 
occurs at site $i$ and that every site $c$ in $J$ is locked at times $(\max(s_{1},s_{3}), s_{2})$. 
Then, the site $i$ is locked at times $(\max(s_{1},s_{3}), s_{2})$, as well.

Our second routine yields a random element of $\A^{ \{a+1,\ldots,b-1\} }$ with
distribution $\prt^{a+1,b-1}(\mu)$. We
make use of a map $\tau : \N \to \N$ such that 
\begin{itemize}
\item $\tau(0)=0$,
\item $\tau(k+1) \ge  \tau(k)+1$ for every $k\ge0$.
\end{itemize}

\begin{tt}
\begin{flushleft}
*** Coalescence detection ***\\
Let $n:=0$; \\
Until $J=\{  a+1,\ldots, b-1 \}$ do: \\
$\{$ Let $J:=\emptyset$; \\

\hspace{.5em} Generate and store $M'_{k}$ and $c'_{k}$ for $k:=\tau(n),\ldots,\tau(n+1)-1$; \\
\hspace{.5em} For $\ell$ going backwards from $\tau(n)$ to $0$ do: \\ 
\hspace{.5em}$\{$ If a coupling event of type $(s_{1}, t'_{\ell} , s_{3},J)$ occurs at site $c'_{\ell}$,
with\\
\hspace{1em}
 $s_{1},s_{3}\le t'_{\tau(n)}$, and if $c'_{\ell}$ is not already in $J$,
then let $J:=J \cup \{ c'_{\ell} \}$;  $\}$ \\ 
\hspace{1em} Let $n:=n+1$; $\}$\\
*** Sampling ***\\
 Let $x := (A,\ldots,A)$; \\
 For $k$  going backwards from $\tau(n)-1$ to $0$ do:  \\
$\{$ Let $x := \gamma(x, M'_{k}, c'_{k})$; \\
\hspace{3ex} Generate $N_{k}$; \\
 \hspace{3ex} For $m$ going from $1$ to $N_{k}$ do: 
\\
\hspace{6ex} $\{$ Generate
 $M''_{k,m}$ and $c''_{k,m}$ ; 
Let $x:=  \gamma(x, M''_{k,m}, c''_{k,m})$; $\}$    $\}$ 
\\
  
Return $\prt^{a+1,b-1}(x)$.
 
 \end{flushleft}
 \end{tt}

As before, the feasibility of this routine relies first on the fact
that the various random quantities can be easily simulated when needed
by the algorithm.  The second key point is that, going forward from
time $-t'_{\tau(n)}$ to time $0$, one can easily detect the coupling
events of type $(s_{1}, t'_{\ell} , s_{3},J)$ with $s_{1}$ and $s_{3}
\le t'_{\tau(n)}$, simply by recording for each site $i$ the latest
(with time going forward) move that affected this site before the
current time $-t'_{\ell}$. The second pass of the routine, that is,
computing $x$ once the coalescence has been obtained, is the same as
in the first method.

Despite the fact that this second routine may use several passes to
detect coalescence, instead of a single one as in the first routine,
the use of sites that are already locked, to detect coupling events,
makes it more effective in some cases, the coalescence time being
always smaller than in the first routine.  The choice of the updating
policy contained in the function $\tau(\cdot)$ is still subject to
empirical adjustment.  The aim here is to reduce the number of passes
as much as possible, while keeping each pass not too time-consuming.

\begin{remark}
In either algorithm, one never generates the ringing times which govern the dynamics themselves, but
only the sequence of moves that they induce.
\end{remark}

\subsection{A special case}
\label{ss.asc}

It is possible to improve upon the previous results in a special case,
namely when  
the only YpR substitutions are from CpG and when the rates of
substitutions are such that 
$$
w_{C}=v_{C},\qquad
w_{G}=v_{G}.
$$
Throughout this section, we assume that these
additional assumptions hold, and we only state the relevant
results, since the proofs closely parallel those in the general case.  

The key observation is the following modification of lemma~\ref{l.d1}, which allows for the use of coarser quotients of the state space $\A$. 
Recall that $\un_{x}(x):=1$ and that $\un_{x}(y):=0$ for every $y\neq x$. 

Fix $\II=\{a,\ldots,b\}$.
\begin{lemma}\label{l.d2}
For every site $i$ in $\II$ and every time $s \ge  \sigma$, 
the functions 
$$\un_{C}\left[\varphi^{\II}_{\sigma} (x,
\xi,\lef(i),s)\right],\quad
\varphi^{\II}_{\sigma} (x, \xi,i,s),\quad
\un_{G}\left[\varphi^{\II}_{\sigma} (x,
\xi,\rig(i),s)\right],
$$
are measurable with respect to
the following initial conditions and source of moves:
$$
\un_{C}(x_{\lef(i)}),\quad
x_{i},\quad
\un_{G}(x_{\rig(i)}),\quad
\xi_{\lef(i)}[\sigma,s],\quad
\xi_{i}[\sigma,s],\quad
\xi_{\rig(i)}[\sigma,s].
$$ 
\end{lemma}
This lemma suggests a modified definition of coupling events.
\begin{definition}[Modified coupling events]\label{d.mce} 
A modified coupling event occurs at site $i$ and times $(s_{1}, s_{2}, s_{3})$ if:
\begin{itemize}
\item $s_{1} > s_{2}$ and $s_{3} > s_{2} $,
\item $-s_{1} $ belongs to $\UU^{A}_{\lef(i)}  \cup \UU^{G}_{\lef(i)} \cup  \UU^{T}_{\lef(i)} \cup \VV^{A}_{\lef(i)}  \cup
\VV^{G}_{\lef(i)}$,
\item $-s_{3} $ belongs to $  \UU^{C}_{\rig(i)}  \cup \UU^{T}_{\rig(i)}  \cup \UU^{A}_{\rig(i)} \cup \VV^{C}_{\rig(i)}  \cup
\VV^{T}_{\rig(i)} $,
\item  $-s_{2} $ belongs to $\bigcup_{z \in \A}  \UU^{z}_{i}$,
\item $(-s_{1},-s_{2}) \cap    \bigcup_{z \in \A }  \UU^{z}_{\lef(i)}\cup\VV^{z}_{\lef(i)} $ is empty,
\item $(-s_{3},-s_{2}) \cap    \bigcup_{z \in \A }  \UU^{z}_{\rig(i)}\cup\VV^{z}_{\rig(i)} $ is empty.
\end{itemize}
\end{definition}
%
The following property is the analogue of proposition~\ref{p:coa}.
\begin{prop}\label{p:coa2}
If a modified coupling event occurs at site $i$ and times  $(s_{1}, s_{2}, s_{3})$, then
$i$ is locked at times  $(s_{4}, s_{2})$, with $s_{4}=\max(s_{1},s_{3})$.   
\end{prop}

\begin{definition}[Modified coupling times]\label{d.mct} 
The modified coupling time $\widetilde T_{i}$ of site $i$ is the supremum of the times $s_{4}$ such that there exists $s_{1}$, $s_{2}$ 
and $s_{3}$, with the following properties:
$s_{2}\ge0$, $s_{4}=\max(s_{1},s_{3})$, and a modified
coupling event occurs at site $i$ and at times $(s_{1}, s_{2}, s_{3})$.
\\
The modified coupling time $\widetilde T_{a,b}$ of the interval
$\II=\{a,\ldots,b\}$
concerns the sites in the interior of $\II$, namely
$$
\widetilde T_{a,b} = \max\{\widetilde T_{i}\,;\,a+1\le i\le b-1\}.
$$
\end{definition}
The definitions below mimick definitions~\ref{d.bata} and \ref{d.nalpha}.
\begin{definition}[Modified combined rates]\label{d.mcr} 
Introduce $$
\vara_{AGT} = \sum_{z = A,G,T} \min(c_{z}, v_{z}), \quad   
\vara_{ACT} = \sum_{z = A,C,T} \min(c_{z}, v_{z}).
$$
and 
$$
\widetilde \bata = \frac{ \vara_{AGT}+ \varb_{R}  }{ \vara + \varb  } \times 
\frac{ \vara_{ACT}+ \varb_{Y}  }{ \vara + \varb  }.
$$
\end{definition}
The arguments in section~\ref{ss.ce} yield the following estimates.
\begin{prop}
For every integer $N$,
$$\Q (\widetilde  T_{a,b}   \ge   N n_{\widetilde\bata} /\vara) \le   3   \left[ 1    -       \left( 1    - 
2\,(1-\widetilde \bata)^{N}  \right)^{k_{a,b}} \right].
$$
For every $t$, 
$$\Q((\widetilde \bata\vara)\, \widetilde T_{a,b}   \ge   t_{a,b}-\log(1-\widetilde \bata)+t ) \le \exp(-t/4).$$
\end{prop}
Since $\widetilde \bata \ge \bata$, the estimates below compare favorably with those obtained in section~\ref{ss.ce}.

To comply with modified coupling events, the two methods of practical detection of the coalescence, 
described in sections~\ref{ss.fa} and \ref{ss.sa}, can be modified in a straightforward way. This yields a priori
shorter coalescence times, but the effective magnitude of this gain should be evaluated, relying on concrete cases.

\section{Index of notations}
\label{s.ion}
\paragraph{Nucleotides}
\begin{description}
\item
$\A=\{A,C,G,T\}$
\item
$A^\ast=G$, $T^\ast=C$, $C^\ast=T$, $G^\ast=A$
\item
$\prp(A)=\prp(G)=R$,
$\prp(C)=\prp(T)=Y$
\item
$\prr(A)=\prr(G)=R$,
$\prr(C)=C$,
$\prr(T)=T$
(the application $\rho$ defined on $\A$ is not to be
confused with the real number $\varrho$, see below)
\item
$\pry(C)=\pry(T)=Y$,
$\pry(A)=A$,
$\pry(G)=G$
\end{description}
\paragraph{Rates of substitutions}
\begin{description}
\item
$v_x$ : rate of the transversions to $x$
\item
$w_x$ : rate of the transition to $x$
\item
$r_x^y$ : rate of the substitution from $yx^\ast$ to $yx$ when $yx$ is a YpR
dinucleotide; rate of the substitution 
from $x^\ast y$ to $xy$ when $xy$ is a YpR
dinucleotide
\item
$\varrho$ : rate of the substitutions $CG\to CA$ and $CG\to TG$ when these
are the only double substitutions and their rates coincide (the real
number $\varrho$ is not to be
confused with the application $\rho$ defined on $\A$, see above)
\end{description}
\paragraph{Functionals of the rates of substitutions}
\begin{description}
\item
$u_x=v_x-w_x$
\item
$v=v_A+v_T+v_C+v_G$,
$w=w_A+w_T+w_C+w_G$
\item
$t_Y=(v_C+v_T)/v$,
$t_R=(v_A+v_G)/v$
\item
$s_A=s_G=s_R=v_A+w_T+w_C+v_G$,
$s_T=s_C=s_Y=w_A+v_T+v_C+w_G$
\item
$v^\ast_A=v_A/s_Y$,
$v^\ast_T=v_T/s_R$,
$v^\ast_C=v_C/s_R$,
$v^\ast_G=v_G/s_Y$
\item
$\vara_R=\min(c_A,v_A)+\min(c_G,v_G)$,
$\vara_Y=\min(c_T,v_T)+\min(c_C,v_C)$,
$\vara=\vara_R+\vara_Y$
\item
$\varb_R=(c_A-v_A)^++(c_G-v_G)^+$,
$\varb_Y=(c_T-v_T)^++(c_C-v_C)^+$,
$\varb=\varb_R+\varb_Y$
\item
$\bata=t_R\,t_Y$,
$n_{\bata}=-4\log(1-\bata)/\bata$
\item
$\vara_{AGT}=\vara_R+\min(c_T,v_T)$,
$\vara_{ACT}=\vara_Y+\min(c_A,v_A)$
\item
$\widetilde\bata=(\vara_{AGT}+\varb_R)\,(\vara_{ACT}+\varb_Y)/(\vara+\varb)^2$
\end{description}
\paragraph{Rates of substitutions in the S/W symmetric case}
\begin{description}
\item
$v_W=v_A=v_T$ when $v_A=v_T$,
$v_S=v_C=v_G$ when $v_C=v_G$
\item
$w_W=w_A=w_T$ when $w_A=w_T$,
$w_S=w_C=w_G$ when $w_C=w_G$
\item
$r_W$ : rate of the substitutions $CG\to CA$ and $CG\to TG$ when these
two rates coincide
\item
$r_S$ : rate of the substitutions $TA\to CA$ and $TA\to TG$ when these
two rates coincide
\item
$\ss_S=v_S+w_S$, $\ss_S=v_W+w_W$
\item
$\sv=v_S+v_W$,
$\sw=w_S+w_W$
\item
$\ss=\sv+\sw=\ss_S+\ss_W$
\end{description}
\paragraph{Sets of sites}
\begin{description}
\item
$\II=\ab$ or $\II=\Z$
\item
$\prt^{\ab}$ : projection on the $\ab$ coordinates
\item
$\lef(i)$ and $\rig(i)$ : left and right neighbors of the site $i$ in $\II$
\item
$\displaystyle k_{a,b}=\pentsup{\frac{b-a-1}3}$,
$t_{a,b}=6\log k_{a,b}$
\end{description}


\addcontentsline{toc}{section}{References}

{\sc Jean B\'erard, Didier Piau}
\begin{quote}
Institut Camille Jordan - UMR 5208
\\ Université Claude Bernard Lyon 1
\\ 43 boulevard du 11 novembre 1918
\\ 69622 Villeurbanne Cedex
 France 
\\
{\tt [Jean.Berard,Didier.Piau]@univ-lyon1.fr}
\\
{\tt lapcs.univ-lyon1.fr/$\sim$[jberard,piau]}
\end{quote}

{\sc Jean-Baptiste Gou\'er\'e}
\begin{quote}
Laboratoire  MAPMO - UMR 6628 
\\ Université d'Orléans 
\\ B.P. 6759
\\ 45067 Orléans Cedex 2
 France
\\ {\tt Jean-Baptiste.Gouere@univ-orleans.fr}
\\ {\tt www.univ-orleans.fr/mapmo/membres/gouere}
\end{quote}

\end{document}